\documentclass{amsart}
\usepackage{epsfig}

\newtheorem{theorem}{Theorem}[section]
\newtheorem{lemma}[theorem]{Lemma}
\newtheorem{prop}[theorem]{Proposition}

\theoremstyle{definition}
\newtheorem{definition}[theorem]{Definition}

\theoremstyle{remark}

\numberwithin{equation}{section}

\def\AA{\mathbb A}

\newcommand{\R}{\mathbb R}

\newcommand{\Z}{\mathbb Z}
\newcommand{\N}{\mathbb N}

\renewcommand{\d}{\delta}
\renewcommand{\l}{\lambda}
\renewcommand{\L}{\Lambda}
\renewcommand{\a}{\alpha}

\newcommand{\D}{\mathcal D}
\newcommand{\cC}{\mathcal C}
\newcommand{\cB}{\mathcal B}
\newcommand{\cS}{\mathcal S}
\newcommand{\Sh}{\mathcal S\! \mathcal H}
\newcommand{\cR}{\mathcal R}

\newcommand{\cQ}{\mathcal Q}

\renewcommand{\H}{\mathcal H}

\newcommand{\W}{\mathcal W}

\newcommand{\supp}{\text{supp }\!}
\newcommand{\sing}{\text{sing }\!}

\newcommand{\ip}[2]{\langle #1 , #2 \rangle}

\newcommand{\norm}[1]{\| #1 \|}

\renewcommand{\qed}{\hskip1em \Box}
\allowdisplaybreaks

\newcommand{\absip}[2]{\left| \left\langle#1,#2\right\rangle \right|}

\begin{document}

\title{Resolution of the Wavefront Set using Continuous Shearlets}

\author{Gitta Kutyniok}
\address{Institute of Mathematics,
Justus--Liebig--University Giessen, 35392 Gies\-sen, Germany}
%\curraddr{Department of Mathematics and Statistics,
%Case Western Reserve University, Cleveland, Ohio 43403}
\email{gitta.kutyniok@math.uni-giessen.de}
\thanks{The first author was supported by DFG research fellowship
KU 1446/5.}

\author{Demetrio Labate}
\address{Department of Mathematics,
North Carolina State University, Campus Box 8205, Raleigh, NC 27695,
USA}
\email{dlabate@unity.ncsu.edu}

%    General info

\date{April 21, 2006.}



\begin{abstract}
It is known that the continuous wavelet transform of a function
$f$ decays very rapidly near the points where $f$ is smooth, while
it decays slowly near the irregular points. This property allows
one to precisely identify the singular support of $f$. However,
the continuous wavelet transform is unable to provide additional
information about the geometry of the singular points. In this
paper, we introduce a new transform for functions and
distributions on $\R^2$, called the Continuous Shearlet Transform.
This is defined by $\Sh_f(a,s,t) = \ip{f}{\psi_{ast}}$, where the
analyzing elements $\psi_{ast}$ are dilated and translated copies
of a single generating  function $\psi$ and, thus, they form an
affine system. The resulting continuous shearlets $\psi_{ast}$ are
smooth functions at continuous scales $a >0$, locations $t \in
\R^2$ and oriented along lines of slope $s \in \R$ in the
frequency domain. The Continuous Shearlet Transform transform  is
able to identify not only the location of the singular points of a
distribution $f$, but also the orientation of their distributed
singularities. As a result, we can use this  transform to
exactly characterize the wavefront set of $f$.
\end{abstract}

\maketitle

\section{Introduction}

It is well-known that if $\psi$ is a `nice' continuous wavelet on
$\R^n$, and $f$ is a function that is smooth apart from a
discontinuity at $x_0 \in \R^n$, then the continuous wavelet
transform
$$ \W_f(a,t) = a^{-\frac{n}{2}} \int_\R f(x) \, \psi\left(a^{-1}(x-t)\right) \, dx, \quad a >0, t \in \R^n$$
decays rapidly as $a \to 0$, unless $t$ is near $x_0$
\cite{Hol95,Mey}. As a consequence, the continuous wavelet transform
is able to resolve the  {\it singular support} of the distribution
$f$, that is, to identify the set of points where $f$ is not
regular. However, the transform $\W_f(a,t)$ is unable to provide
additional information about the {\it geometry} of the set of
singularities of a more general distribution $f$. For example, in
many situations one would like to identify not only the location of
a certain distributed singularity, but also its orientation. This is
very relevant, in particular, in the study of the propagation of
singularities associated with partial differential equations
\cite{Hor,Sog}.

In this paper, we introduce a new two-dimensional continuous
transform, called {\it continuous shearlet transform}, that is
mapping a tempered distribution $f \in \cS'(\R^2)$ to
$\Sh_f(a,s,t)$, where $a>0$,  $s \in \R$ and $t \in \R^2$, are the
scale, shear and location variables, respectively. The transform
is defined by $\Sh_f(a,s,t) = \ip{f}{\psi_{ast}}$, where the
analyzing elements $\psi_{ast}$, called {\it continuous
shearlets}, are dilated and translated copies of a single
generating function $\psi$, thereby constituting a continuous
2-dimensional wavelet-like system. The generator $\psi \in
L^2(\R^2)$ is chosen to be arbitrarily smooth and has compact
support in the frequency domain with a needlelike structure to
capture directionality. For the dilation matrix we employ the
product of a parabolic scaling matrix associated with $a>0$ and a
shear matrix associated with $s \in \R$, whereas $t \in \R^2$
serves as the translation parameter. The system $\{\psi_{ast} : a
> 0,\, s \in \R,\, t \in \R^2\}$ will be shown to satisfy a
Calder\`on-type formula with respect to a special measure, that
is, it forms a reproducing system for $L^2(\R^2)$.

We will prove that the continuous shearlet transform is able  to
identify not only the singular support of the distribution $f$,
but also the orientation of distributed singularities along
curves. In particular, the decay properties of the continuous
shearlet transform as $a \to 0$  precisely characterize the {\it
wavefront set} of $f$ (see Section~\ref{s.wav}).

Historically, the idea of using wavelet-like transforms to
identify the set of singularities of a distribution can be traced
back to the {\it wave packet transform}, introduced independently
by Bros and Iagolnitzer \cite{BI76} and C\'{o}rdoba and Fefferman
\cite{CF78}. More recently, Smith \cite{Smi98} and Cand\`es and
Donoho \cite{CD04a,CD04b} have introduced a continuous transform
that uses parabolic scaling in polar coordinates and has the
ability to identify the wavefront set of a distribution. In
particular, the Continuous Curvelet Transform (CCT) of Cand\`es
and Donoho has properties similar to the continuous shearlet
transform presented in this paper. Namely, in both cases, the
wavefront set is exactly characterized by the decay properties of
the continuous transform. Unlike the CCT, however, our approach
exploits the framework of {\it affine systems}. As a consequence,
the continuous shearlet transform  is not only much closer in
spirit to the traditional continuous wavelet transform, but also
avoids the more complicated structure of the curvelet
construction, which uses infinitely many different generators.

Another motivation for this investigation and the use of the
framework of affine systems comes from the study of discrete
wavelets, and, more specifically their ability to approximate
efficiently smooth functions with singularities. This property is
closely related to the micro-localization property of the continuous
wavelet transform. To illustrate this point, let $f$ be a
one-dimensional function  that is smooth apart from a discontinuity at
$x_0$
 and consider the wavelet representation of $f$:
 $$ f = \sum_{k, j \in \Z} \ip{f}{\psi_{j,k}} \, \psi_{j,k}(x), $$
where  $\psi_{j,k}(x) = 2^{-j/2} \, \psi(2^{-j} x-k)$ and $\psi$ is
a `nice' wavelet.
 Since the wavelet transform $\W_f(2^{-j},k) = \ip{f}{\psi_{j,k}}$ decays rapidly for $j \to \infty$ unless $k$ is near $x_0$,
it follows that very few coefficients of the wavelet representation
are sufficient to approximate $f$ accurately. Indeed, the wavelet
representation is optimally sparse for this type of functions
(cf.~\cite[Ch.9]{Mal98}). Therefore, the properties of the
continuous wavelet transform are crucial in order to achieve a
sparse representation.

Traditional wavelets, however, have a very limited capability in
dealing with discontinuities in higher dimensions. Consider, for
example, the wavelet representation of a 2-dimensional function
that is smooth away from a discontinuity along a  curve. Because
the discontinuity is spatially distributed, it interacts
extensively with the elements of the wavelet basis, and, thus,
many wavelet coefficients are needed to represent $f$ accurately.
Related to this is the fact that the corresponding two-dimensional
continuous wavelet transform is unable to `track' the
discontinuous curve. As pointed out by several authors
(cf.~\cite{CD99,CD04}), in order to represent multidimensional
functions efficiently, one has to use representations that are
much more flexible than traditional wavelets, and have the ability
to capture the geometry of multidimensional phenomena. One goal of
this paper is to show that the study of the continuous wavelet
transform associated with the affine group on $\R^2$ provides a
unifying perspective on the construction of such representations.

The study of the discrete analog of the continuous shearlet
transform is currently being developed by the authors and their
collaborators \cite{GKL05,GL,GLLWW,GLLWW3,LLKW05}. In particular,
it was shown that, thanks to the mathematical structure of affine
systems, discrete shearlets are associated with a Multiresolution
analysis similar to traditional wavelets, which is very relevant
for the development of fast algorithmic implementations. In
addition, shearlets provide optimally sparse representations for
functions in certain classes and can be easily generalized to
higher dimensions.

\medskip
The paper is organized as follows. In Section~\ref{s.adm} we
recall the basic properties of affine systems on $\R^n$ and the
continuous wavelet transform, before introducing the continuous
shearlet transform  in Section~\ref{s.adm}. In Section~\ref{s.cst}
we apply this new transform to several examples of distributions
containing different types of singularities. The main result of
this paper is proved in Section~\ref{s.wav}, where we show that
the continuous shearlet transform exactly characterizes the
wavefront set of a distribution. Finally, in Section~\ref{s.ext},
we discuss several variants and generalizations of our
construction.

\section{Affine systems and wavelets}  \label{s.adm}

\subsection{One-dimensional Continuous Wavelet Transform}

Let $\AA_1$ be the {\it affine group} associated with $\R$,
consisting of all pairs $(a,t)$, $a,t  \in \R, a > 0$, with group
operation $ (a,t) \cdot (a',t')   = (a a', t + a t')$. The {\it
(continuous) affine systems} generated by $\psi \in L^2(\R)$ are
obtained from the action of the {\it quasi--regular representation}
$\pi_{(a,t)}$ of $\AA_1$ on $L^2(\R)$, that is
$$ \bigr\{ \psi_{a,t}(x) = \pi_{(a,t)} \, \psi (x) = T_t \, D_a \,  \psi(x): \, (a,t) \in \AA_1 \bigl\},$$
where the {\it translation operator} $T_t$ is defined by
$T_{t}\psi(x)=\psi(x-t)$ and the {\it  dilation operator} $D_a$ is
defined by $D_a\psi(x)=a^{-1/2} \psi(a^{-1} x)$.

It was observed by Calder\`on \cite{Cal64} that, if $\psi$ satisfies
the {\it admissibility} condition
\begin{equation}\label{eq.ac}
\int_{0}^\infty |\hat \psi (a \xi)|^2 \, \frac{ da }{a} = 1 \quad
\text{ for a.e. } \xi \in \R,
\end{equation}
then any $f \in L^2(\R)$ can be recovered via the reproducing
formula:
$$ f = \int_{\AA_1} \ip{f}{\psi_{a,t}} \, \psi_{a,t} \,  d \mu(a,t), $$
where $d \mu(a,t) = d t \, \tfrac{da}{a^2} $ is the left Haar
measure of $\AA_1$. Here the Fourier transform is defined by
$\hat{\psi}(\xi) = \int \psi(x) \, e^{-2 \pi i \xi  x} \, d x$. As
usual, $\check{\psi}$ will denote the inverse Fourier transform. The
function $\psi$ is called a  {\it continuous wavelet}, if $\psi$
satisfies \eqref{eq.ac}, and $\W_f(a,t) = \ip{f}{\psi_{a,t}}$ is the
{\it continuous wavelet transform} of $f$. We refer to~\cite{GMP85}
for more details about this.

Discrete affine systems and wavelets are obtained by
`discretizing' appropriately the corresponding continuous systems.
In fact, by replacing  $(a,t) \in \AA_1$ with the discrete set
$(2^j, 2^j m)$, $j,m \in \Z$, one obtains the discrete dyadic affine
system
\begin{equation}\label{eq.ddas}
\bigr\{ \psi_{j,m}(x) = T_{2^j m} \, D_2^j \, \psi(x) =  D_2^j \,
T_m \, \psi(x): \, j,m \in \Z \bigl\},
\end{equation}
and $\psi$ is called a {\it wavelet} if (\ref{eq.ddas}) is an
orthonormal basis or, more generally, a Parseval frame for
$L^2(\R)$.

Recall that a countable collection $\{\psi_i\}_{i \in I}$ in a
Hilbert space $\H$ is a {\it Parseval frame} (sometimes called a
{\it tight frame}) for $\H$, if $\sum_{i \in I} |\ip{f}{\psi_{i}}|^2
= \norm{f}^2$ for all $f \in \H$. This is equivalent to the
reproducing formula $f=\sum_{i \in I} \ip{f}{\psi_{i}}\,\psi_{i}$
for all $f \in \H$, where the series converges in the norm of $\H$.
Thus Parseval frames provide basis-like representations even though
a Parseval frame need not be a basis in general. We refer the reader
to~~\cite{Cas00,Chr03} for more details about frames.

%***************************************************************************************************

\subsection{Higher-dimensional Continuous Wavelet Transform}

The natural way of extending the theory of affine systems to higher
dimensions is by replacing $\AA_1$ with the {\it full affine group
of motions on $\R^n$}, $\AA_n$, consisting of the pairs $(M,t) \in
GL_n(\R) \times \R^n$ with group operation $(M,t)\cdot (M',t') = (M
M', t+ M t')$. Similarly to the one-dimensional case, the affine
systems generated by $\psi \in L^2(\R^n)$ are given by
$$ \bigr\{ \psi_{M,t}(x) =  T_t \, D_M \,  \psi(x): \, (M,t) \in \AA_n \bigl\},$$
where here the {\it  dilation operator} $D_M$ is defined by $D_M \,
\psi(x)=|\det M|^{-\frac12} \psi(M^{-1} x)$. The generalization of
the Calder\`on admissibility condition to higher dimensions and the
construction of multidimensional wavelets is a far more complex task
than the corresponding one-dimensional problem, and yet not fully
understood. We refer to~\cite{LWWW,WW} for more details.

Now let $G$ be a subset  of $GL_n(\R)$ and define $\L \subseteq
\AA_n$ by $\L = \{(M,t) : M \in G, t \in \R^n\}$. If there exists a
function $\psi \in L^2(\R^n)$ such that, for all $f \in L^2(\R^n)$,
we have:
\begin{equation}\label{eq.calderon}
f = \int_{\R^n} \int_{G} \ip{f}{T_t \, D_M \,  \psi} \, T_t \, D_M
\, \psi \,  d \l(M) \, d t,
\end{equation}
where  $\l$ is a measure on $G$, then $\psi$ is a {\it continuous
wavelet} with respect to $\L$. The following result, that is a
simple modification of Theorem~2.1 in \cite{WW}, gives an exact
characterization of all those $\psi \in L^2(\R^n)$ that are
continuous wavelets with respect to $\L$. The proof of this theorem is reported in the Appendix.

\begin{theorem} \label{th.rf}
Equality (\ref{eq.calderon}) is valid for all $f \in L^2(\R^n)$ if
and only if, for all $\xi \in {\R}^n \setminus \{0\}$,
\begin{equation}\label{eq.admis}
 \Delta(\psi)(\xi) = \int_G |\hat \psi(M^t \xi) |^2 \, |\det M| \, d \l(M) = 1.
\end{equation}
\end{theorem}
The choice of the measure $\l$ on $G$ is not unique. If $G$ is not
simply a subset of $GL_n(\R)$, but also a subgroup, then we can use
the left Haar measure on $G$ which is unique up to a multiplicative
constant. Also, observe that Theorem~\ref{th.rf} extends to
functions on subspaces of $L^2(\R^n)$ of the form
$$L^2(V)^\vee= \{f \in L^2(\R^n): \supp \hat f \subset V\}.$$

%************************************************************************************************

\subsection{Localization of Wavelets}

The decay properties of the functions $\psi_{M,t} = T_t \, D_M \,
\psi$, where $\hat \psi \in C_0^\infty$, are described by the
following proposition.

\begin{prop} \label{pro.loc}
Suppose that $\psi \in L^2(\R^n)$ is such that $\hat \psi \in
C^\infty_0(R)$, where $R = \supp \hat \psi\subset \widehat \R^n$.
Then, for each $k \in \N$, there is a constant $C_k$ such that, for
any $x \in \R^n$, we have
\begin{equation*}
    |\psi_{M,t}(x)| \le C_k \, |\det M|^{-\frac12} \, (1+| M^{-1} (x-t)|^2)^{-k}.
\end{equation*}
In particular, $C_k = k \, m(R)  \, \bigl(\norm{\hat \psi}_\infty +
\norm{\triangle^k \hat \psi}_\infty \bigr)$, where $\triangle =
\sum_{i=1}^n \frac{\partial^2}{\partial \xi_i^2}$ is the frequency
domain Laplacian operator and $m(R)$ is the Lebesgue measure of $R$.
\end{prop}

The proof of this proposition relies on the following known
observation, whose proof is included for completeness.

\begin{lemma} \label{lem.loc}
Let $g$ be such that $\hat g \in C_0^\infty(R)$, where $R \subset
\R^n$ is the $\supp \hat g$. Then, for each $k \in \N$, there is a
constant $C_k$ such that for any $x \in \R^n$
\begin{equation*}
    |g(x)| \le C_k \, (1+|x|^2)^{-k}.
\end{equation*}
In particular, $C_k = k \, m(R) \, \bigl(\norm{\hat g}_\infty +
\norm{\triangle^k \hat g}_\infty \bigr)$.
\end{lemma}

\noindent {\bf Proof.} Since $g(x) = \int_{R} \hat g (\xi) \, e^{2
\pi i \xi x} \, d \xi$, then, for every $x \in \R^2$,
\begin{equation} \label{e1-lem1}
|g(x)| \le m(R) \, \norm{\hat g}_\infty.
\end{equation}
An integration by parts shows that
$$\int_{R} \triangle \hat g (\xi) \, e^{2 \pi i \xi x} \, d \xi = - (2 \pi)^2 \, |x|^2 \, g(x)$$
and thus, for every $x \in \R^2$,
\begin{equation} \label{e2-lem1}
(2 \pi \, |x|)^{2k} \, |g(x)| \le m(R) \, \norm{\triangle^k \hat
g}_\infty.
\end{equation}
Using (\ref{e1-lem1}) and (\ref{e2-lem1}), we have
\begin{equation} \label{e3-lem1}
\bigl(1+( 2 \pi \, |x|)^{2k} \bigr) \, |g(x)| \le m(R) \, \bigl(
\norm{\hat g}_\infty +\norm{\triangle^k \hat g}_\infty \bigr).
\end{equation}
Observe that, for each $k \in \N$,
$$ (1+|x|^2)^{k} \le \bigl(1+( 2 \pi)^2 \, |x|^2\bigr)^k \le k \, \bigl(1+( 2 \pi \, |x|)^{2k} \bigr). $$
Using this last inequality and (\ref{e3-lem1}), we have that for
each $x \in \R^n$
\begin{equation*}
 |g(x)| \le k \, m(R) \, (1+|x|^2)^{-k}  \, \bigl( \norm{\hat g}_\infty +\norm{\triangle^k \hat g}_\infty \bigr). \qed
\end{equation*}

\medskip
A simple re-scaling argument now proves Proposition~\ref{pro.loc}.

\medskip
{\bf Proof of  Proposition~\ref{pro.loc}.} A direct computation
gives:
\begin{eqnarray*}
\psi( M^{-1}(x-t)) & = & \int_{R} \hat \psi (\xi) \, e^{2 \pi i  M^{-1}(x-t) \xi } \, d \xi \\
& = & \int_{R} \hat \psi (\xi) \, e^{2 \pi i (x-t) M^{-t} \xi  } \, d \xi \\
& = & \int_{R  M^{-1}} \hat \psi (M^t \eta ) \, e^{2 \pi i (x-t)
\eta} \, |\det M| \,  d \eta.
\end{eqnarray*}
It follows that
$$ |\psi( M^{-1}(x-t))| \le m(M^{-1} R  ) \, |\det M| \,  \norm{\hat \psi (M^t \cdot  )}_\infty
= m(R) \, \norm{\hat \psi }_\infty.$$ Using a simple modification of
the argument  in Lemma~\ref{lem.loc}, we have that
$$ (2 \pi \, | M^{-1}(x-t)|)^{2k} \, |\psi( M^{-1}(x-t))| \le m(R) \, \norm{\triangle^k \hat \psi}_\infty.$$
Next, arguing again as in Lemma~\ref{lem.loc} we have  that
$$  |\psi( M^{-1}(x-t))| \le k \, m(R) \, (1+| M^{-1}(x-t)|^2)^{-k}  \, \bigl( \norm{\hat \psi}_\infty +\norm{\triangle^k \hat \psi}_\infty \bigr).$$
This completes the proof. $\qed$
%The proof now follows from the fact that $\psi_{A,t}(x) = |\det A|^{-1/2}  \, \psi( A^{-1}(x-t))$.

%*******************************************************************************************************

\section{Continuous Shearlet Transform} \label{s.cst}

\subsection{Definition}

In this paper, we will be  interested in the affine systems obtained
when $\L$ is a subset of $\AA_2$ of the form
\begin{equation}\label{def.G}
\L = \{ (M,t): \, M \in G, t \in \R^2 \},
\end{equation}
and $G \subset GL_2(\R)$ is the set of matrices:
\begin{equation}\label{def.G2}
G = \left\{ M = M_{as} = \begin{pmatrix}  a & - \sqrt a \, s \\
\\ 0  & \sqrt a  \end{pmatrix}, \quad a \in I, \, s \in S \right\},
\end{equation}
where $I \subset \R^+$, $S \subset \R$.   It is  useful to notice
that the matrices $M$ can be factorized as $M= B \, A$, where $B$ is
the {\it shear matrix}
$ B =  \left(\begin{smallmatrix}  1 & -s \\
\\ 0  & 1  \end{smallmatrix}\right)$ and $A$ is the
diagonal matrix $ A = \left(\begin{smallmatrix}  a & 0 \\
\\ 0  & \sqrt a  \end{smallmatrix}\right).$ In particular, $A$
produces {\it parabolic scaling}, that is, $ f(A x) = f \left(A \left(\begin{smallmatrix}  x_1 \\
\\ x_2  \end{smallmatrix}\right) \right)$ leaves invariant the parabola $x_1= x_2^2$.
 Thus, the matrix $M$ can be interpreted as the superposition of
 parabolic scaling and shear transformation.

We will now construct a continuous wavelet on $\R^2$ associated with
the subset $\L$ of the affine group, given by (\ref{def.G}). We will
consider two situations, corresponding to $I = \R^+$, $S = \R$ or $I
= \{a: 0 \le a \le 1\}$, $S= \{s\in \R: |s| \le s_0\}$, for some
$s_0>0$.

For $\xi = (\xi_1, \xi_2) \in \widehat \R^2$, $\xi_2 \ne 0$, let
$\psi$ be given by
 \begin{equation}\label{ex1_psi}
 \hat \psi (\xi) = \hat \psi (\xi_1,\xi_2) = \hat \psi_1(\xi_1) \,  \hat
 \psi_2(\tfrac{\xi_2}{\xi_1}).
\end{equation}

\begin{prop} \label{prop1}
Let $\L$ be given by (\ref{def.G}) and (\ref{def.G2}) with $I=\R^+$,
$S = \R$, and $\psi \in L^2(\R^2)$ be given by (\ref{ex1_psi})
where:
\begin{enumerate}
    \item[(i)] $\psi_1 \in L^2(\R)$  satisfies
the Calder\`on condition (\ref{eq.ac});
    \item[(ii)] $\norm{\psi_2}_{L^2}=1$.
\end{enumerate}
 Then $\psi$ is a continuous wavelet for $L^2(\R^2)$ with respect to $\L$.
\end{prop}

 {\bf Proof.} A direct computation shows that $M^t (\xi_1,\xi_2)^t=
(a \xi_1,a^{1/2} (\xi_2-s \xi_1  ))^t$. By choosing as measure $d
\l(M) = \frac{d a}{|\det M|^2} \, d s$,
 the admissibility  condition (\ref{eq.admis}) becomes
 \begin{equation}\label{eq.admiA}
 \Delta(\psi)(\xi) = \int_\R \int_{\R^+}
 |\hat \psi_1(a \, \xi_1)|^2 \, |\hat \psi_2(a^{-\frac12} (\tfrac{\xi_2}{\xi_1}-s)) |^2 \, \, a^{-\frac32} d a \, d s = 1.
\end{equation}
Thus, by Theorem~\ref{th.rf}, to show that $\psi$ is a continuous
wavelet it is sufficient to show that  (\ref{eq.admiA}) is
satisfied. Using the assumption on $\psi_1$ and $\psi_2$, we have:
\begin{eqnarray*}
 \Delta(\psi)(\xi) & = & \int_\R \int_{\R^+}
 |\hat \psi_1(a \, \xi_1)|^2 \, |\hat \psi_2(a^{-\frac12} (\tfrac{\xi_2}{\xi_1}-s)) |^2 \, \, a^{-\frac32} d a \, d s  \\
  &=  & \int_{\R^+} |\hat \psi_1(a \, \xi_1)|^2 \, \Bigl( \int_{\R } |\hat \psi_2(a^{-\frac12} \, \tfrac{\xi_2}{\xi_1} -s)|^2
   \, \, d s \Bigr) \, \frac{d  a}{a}\\
  & = & \int_{\R^+} |\hat \psi_1(a\, \xi_1)|^2 \, \frac{d  a}{a}  = 1 \quad \text{ for a.e. } \xi = (\xi_1, \xi_2) \in \R^2.
  \end{eqnarray*}
This shows that equality (\ref{eq.admiA}) is satisfied and, hence,
$\psi$ is a continuous wavelet. $\qed$

\medskip
If the set $S$ is not all of $\R$, then we need some additional
assumptions on $\psi$. Consider the subspace of $L^2(\R^2)$ given by
$L^2(C)^\vee = \{f \in L^2(\R^2): \supp \hat f \subset C\}$, where
\begin{equation*}
C = \{(\xi_1, \xi_2) \in \R^2: |\xi_1| \ge 2 \text{ and }
|\tfrac{\xi_2}{\xi_1}| \le 1\}.
\end{equation*}
We have the following result.

\begin{prop} \label{prop2}
Let $\L$ be given by (\ref{def.G}) and (\ref{def.G2}) with $I = \{a:
0 \le a \le 1\}$, $S = \{s\in \R: |s|\le 2\}$, and $\psi \in
L^2(\R^2)$ be given by (\ref{ex1_psi}) where:
\begin{enumerate}
    \item[(i)]  $\psi_1 \in L^2(\R)$ satisfies
the Calder\`on condition (\ref{eq.ac}), and $\supp \hat \psi_1
\subset [-2,-\frac12] \cup [\frac12,2]$;
    \item[(ii)]  $\norm{\psi_2}_{L^2}=1$ and $\supp \hat
\psi_2 \subset [-1,1]$.
\end{enumerate}
 Then $\psi$ is a continuous wavelet for $L^2(C)^\vee$ with respect to $\L$, that is,  for all $f \in L^2(C)^\vee$,
\begin{equation*}
f(x) = \int_{\R^2} \int_{-2}^2 \int_{0}^1 \ip{f}{ \psi_{a s t}}
\,\psi_{a s t}(x)  \, \frac{d a}{a^3} \, ds \, dt.
\end{equation*}
\end{prop}

{\bf Proof.} We apply again Theorem~\ref{th.rf}, for functions on
$L^2(C)^\vee$. Using the assumptions on $\psi_2$, $S$ and $I$ we
have that, for $\xi \in C$:
$$ \int_{\frac 1{\sqrt{a}}(\frac{\xi_2}{\xi_1}-2)}^{\frac 1{\sqrt{a}}(\frac{\xi_2}{\xi_1}+2)}  |\hat \psi_2 (s)|^2 \, d
s= \int_{-1}^1 |\hat \psi_2(s)|^2 \, d s =1.$$ Thus, for a.e. $\xi
\in C$ we have that
\begin{eqnarray*}
\Delta(\psi)(\xi) & = &  \int_{-2}^2 \int_0^1  |\hat \psi (M_{as}^t \xi)|^2 \, a^{-\frac32} \, d a \, d s \\
& = &   \int_{-2}^2 \int_0^1 |\hat \psi_1 (a \xi_1)|^2 \, |\hat
\psi_2 (a^{-\frac12}(\tfrac{\xi_2}{\xi_1}-s))|^2
a^{-\frac32} \, d a \, d s \\
& = &   \int_0^1 |\hat \psi_1 (a \xi_1)|^2 \, \int_{\frac
1{\sqrt{a}}(\frac{\xi_2}{\xi_1}-2)}^{\frac 1{\sqrt{a}}
(\frac{\xi_2}{\xi_1}+2)}  |\hat \psi_2 (s)|^2 \, d s  \, \frac{d a}{a}  \\
& =  & \int_0^1 |\hat \psi_1 (a \xi_1)|^2  \, \frac{d a}{a}.
\end{eqnarray*}
Since $\xi_1 \ge 2$, using the assumptions on the support of $\hat
\psi_1$ and condition~(\ref{eq.ac}), from the last expression we
have that, for a.e. $\xi \in C$,
$$ \Delta(\psi)(\xi) =  \int_0^{\xi_1} |\hat \psi_1 (a)|^2  \, \frac{d a}{a} = \int_{\frac 12}^{2} |\hat \psi_1 (a)|^2  \, \frac{d
a}{a} = \int_{0}^{\infty} |\hat \psi_1 (a)|^2  \, \frac{d a}{a}
=1.$$ This shows that the admissibility condition (\ref{eq.admis})
for this system is satisfied and this completes the proof. $\qed$

\medskip
There are several examples of functions $\psi_1$ and $\psi_2$
satisfying the assumptions of Proposition~\ref{prop1} as well as
Proposition~\ref{prop2}. In addition, we can choose $\psi_1, \psi_2$
such that
 $\hat \psi_1$, $\hat \psi_2 \in C_0^\infty$ (see~\cite{GKL05,GLLWW3} for the construction of these functions).

%\vspace*{-0.5cm}

\begin{figure}[h]
\begin{center}
\begin{picture}(300,150)(0,0)
  \put(80,0){\epsfig{file=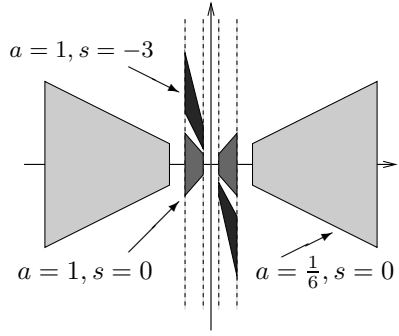,width=5cm}}
  \put(120,30){\vector(1,1){19}}
  \put(78,20){\text{$a=1,s=0$}}
  \put(120,100){\vector(2,-1){19}}
  \put(75,104){\text{\small $a=1,s=-3$}}
  \put(178,28){\vector(2,1){19}}
  \put(168,20){\text{$a=\frac16, s=0$}}
  \end{picture}
\end{center}
\caption{Support of the shearlets $\hat{\psi}_{a s t}$ (in the
frequency domain) for different values of $a$ and $s$.}
\end{figure}

\medskip
Now we can define the {\it continuous shearlet transform:}

\begin{definition}
\label{defi.shear}  Let $\psi \in L^2(\R^2)$ be given by
(\ref{ex1_psi}) where:
\begin{enumerate}
    \item[(i)]  $\psi_1 \in L^2(\R)$ satisfies the Calder\`on condition (\ref{eq.ac}), and
    $ \hat \psi_1 \in C^\infty_0(\R)$ with $\supp \hat \psi_1 \subset [-2,-\frac12] \cup [\frac12,2]$;
    \item[(ii)] $\norm{\psi_2}_{L^2}=1$, and  $\hat \psi_2 \in C^\infty_0(\R)$ with $\supp \hat
\psi_2 \subset [-1,1]$ and $\hat \psi_2 > 0$ on $(-1,1)$.
\end{enumerate}
The functions generated by $\psi$ under the action of $\L$, namely:
\[    \psi_{ast}(x) =  T_t \, D_M \psi(x) = a^{-\frac34} \, \psi \left( \left( \begin{smallmatrix}  a & -\sqrt a \, s \\
\\ 0  & \sqrt a  \end{smallmatrix}\right)^{-1}  (x-t)\right),\]
where $a \in I \subset \R^+$,  $s \in S \subset \R$ and $t \in \R^2$,
are called {\em continuous shearlets}. The {\em continuous shearlet
transform} of $f$ is defined by
$$ \Sh_f(a,s,t) = \ip{f}{\psi_{ast}}, \, \quad a\in I \subset \R^+, s \in S \subset \R, t \in \R^2. $$
\end{definition}

Observe that, unlike the traditional wavelet transform  which depends
only on scale and translation, the shearlet transform is a function
of three variables, that is, the scale $a$, the shear $s$ and the
translation $t$. Many properties of the continuous shearlets are
more evident in the frequency domain. A direct computation shows
that
\begin{eqnarray*}
\hat{\psi}_{ast}(\xi)
 & = &  a^{\frac34} \, e^{-2 \pi i \xi t} \, \hat \psi ( a \, \xi_1,  \sqrt a (\xi_2-s \, \xi_1 )) \nonumber \\
& = &  a^{\frac34} \, e^{-2 \pi i \xi t} \, \hat \psi_1 ( a \,
\xi_1) \, \hat \psi_2  (a^{-\frac12} (\tfrac{\xi_2}{\xi_1}-s)).
\end{eqnarray*}
Thus, each function $\hat{\psi}_{ast}$ is supported in the set:
$$ \supp \hat{\psi}_{ast} \subset \{(\xi_1, \xi_2):
\xi_1 \in [-\tfrac{2}{a}, -\tfrac{1}{2a}] \cup [\tfrac{1}{2a},
\tfrac{2}{a}], \, |\tfrac{\xi_2}{\xi_1}-s| \le \sqrt a\}.$$ As
illustrated in Figure~1, each continuous shearlet ${\psi}_{ast}$ has
frequency support on a pair of trapezoids, symmetric with respect to
the origin, oriented along a line of slope $s$. The support becomes
increasingly thin as $a \to 0$.

\medskip

When $S =\R$ and $I= \R^+$, by Proposition~\ref{prop1}, the
continuous shearlet transform provides a reproducing
formula~(\ref{eq.calderon}) for all $f \in L^2(\R^2)$:
$$ \norm{f}^2 = \int_{\R^2} \int_{-\infty}^\infty \int_{0}^\infty \left| \Sh_f(a,s,t)\right|^2 \,  \frac{d a}{a^3} \, ds \, dt.$$
On the other hand, if $S$, $I$ are bounded sets, by
Proposition~\ref{prop2}, the continuous shearlet transform provides
a reproducing formula only for functions in a proper subspace of
$L^2(\R^2)$. However, even when $S, I$ are bounded, it is possible
to obtain a reproducing formula for all $f \in L^2(\R^2)$ as
follows. Let
$$  \hat \psi^{(v)} (\xi) = \hat \psi^{(v)}  (\xi_1,\xi_2) = \hat \psi_1(\xi_2) \,  \hat
 \psi_2(\tfrac{\xi_1}{\xi_2}),$$
where $\hat \psi_1$, $\hat \psi_2$ are defined as in Definition
\ref{defi.shear}, and let $\L^{(v)} = \{(M,t): M \in G^{(v)}, t
\in \R^2\}$, where
\begin{equation}\label{def.G2v}
G^{(v)} = \left\{ M = M_{as} = \begin{pmatrix}  \sqrt a & 0  \\
\\ - \sqrt a \, s  & a  \end{pmatrix}, \quad a \in I, \, s \in S \right\},
\end{equation}
 Then, proceeding as above, it is easy to show that $\psi^{(v)}$ is a
continuous wavelet for $L^2(C^{(v)})^\vee$ with respect to
$\L^{(v)}$, where $C^{(v)}$ is the vertical cone:
\begin{equation*}
C^{(v)} = \{(\xi_1, \xi_2) \in \R^2: |\xi_2| \ge 2 \text{ and }
|\tfrac{\xi_2}{\xi_1}| > 1\}.
\end{equation*}
Accordingly, we introduce the shearlets $\psi^{(v)}_{ast} = T_t \,
D_{M} \psi^{(v)}$, for $(M,t) \in \L^{(v)}$, and the associated
continuous shearlet transform $ \Sh_f^{(v)}(a,s,t) =
\ip{f}{\psi^{(v)}_{ast}}$. Finally, let $W(x)$ be such that $\hat
W(\xi) \in C^\infty(\R^2)$ and
\begin{equation}\label{eq.W}
    |\hat W(\xi)|^2 + \chi_{C_1}(\xi)  \int_0^1 |\hat \psi_1(a \xi_1)|^2 \, \frac{d a}{a} + \chi_{C_2}(\xi)  \int_0^1 |\hat \psi_1(a \xi_2)|^2 \, \frac{d a}{a} = 1,
\end{equation}
 for a. e. $\xi \in \R^2$, where $C_1 = \{(\xi_1, \xi_2) \in
\R^2: |\tfrac{\xi_2}{\xi_1}| \le 1\}$, $C_2 = \{(\xi_1, \xi_2) \in
\R^2: |\tfrac{\xi_2}{\xi_1}| > 1\}$. Then it follows that $W$ is a
$C^\infty$ window function in $\R^2$ with $\hat W(\xi) = 1$ for
$\xi \in [-1/2,1/2]^2$, $\hat W(\xi) = 0$ outside the box $\{\xi
\in [-2,2]^2\}$. Finally, let $(P_{C_1} f)^\wedge = \hat f \,
\chi_{C_1}$ and $(P_{C_2} f)^\wedge = \hat f \, \chi_{C_2}.$ Then,
for each $f \in L^2(\R^2)$ we have:
\begin{eqnarray}
\norm{f}^2 &=& \int_{\R} |\ip{f}{T_t \, W}|^2 \, d t + \int_{\R^2}
\!\int_{-2}^2 \int_{0}^1
|\Sh_{(P_{C_1} f)}(a,s,t)|^2 \, \frac{d a}{a^3} \, ds \, dt \nonumber \\
&+& \int_{\R^2} \! \int_{-2}^2 \int_{0}^1
|\Sh_{(P_{C_2}f)}^{(v)}(a,s,t)|^2 \, \frac{d a}{a^3} \, ds \, dt. \label{eq.dec3}
\end{eqnarray}
The proof of this equality is reported in the Appendix.
Equation (\ref{eq.dec3}) shows that
$f$ is continuously reproduced  by using isotropic window
functions at coarse scales, and two sets of continuous shearlets at
fine scales: one set corresponding to the horizontal cone $C$ (in
the frequency domain) and another set corresponding to the vertical
cone $C^{(v)}$. The advantage of this construction, with respect to
the simpler one where $S=\R$, is that in this case the set $S$
associated with the shear variable is the closed interval $S = \{s:
|s| \le 2\}$. This property will be important in
Subsection~\ref{s.curv} and Section~\ref{s.wav}.

\bigskip
There are other choices of the subset $\L$, given by~(\ref{def.G}),
generating affine systems with properties similar to the continuous
shearlets. Variants and generalizations of this construction  will
be discussed in Section~\ref{s.ext}.

%***************************************************************************************************

\subsection{Localization of Shearlets}

Since the continuous wavelet $\psi$ associated with the continuous
shearlet transform satisfies $\hat \psi \in C_0^\infty(\widehat
\R^2)$, it follows that the continuous shearlets  {\it decay rapidly
} as $|x| \to \infty$, that is:
$$\psi_{ast}(x)  = O(|x|^{-k}) \quad \text{ as } |x| \rightarrow \infty, \quad \text{ for every } k \ge 0.$$

More precisely, we have the following result.

\begin{prop} \label{pro.locshear}
Let $\psi \in L^2(\R^2)$ be the continuous wavelet associated with the
continuous shearlet transform, and let $M$ be defined as in
\eqref{def.G2}. Then, for each $k \in \N$, there is a constant $C_k$
such that, for any $x \in \R^2$, we have
\begin{eqnarray*}
|\psi_{ast}(x)| & \le & C_k \, |\det M|^{-\frac12} \, (1+| M^{-1}
(x-t)|^2)^{-k}\\
& = & C_k \, a^{-\frac34} \,
(1+a^{-2}(x_1-t_1)^2+2a^{-2}s(x_1-t_1)(x_2-t_2)\\
&&    +a^{-1}(1+a^{-1}s^2)(x_2-t_2)^2)^{-k}.
\end{eqnarray*}
In particular, $C_k = k \, \frac{15}{2} \frac{\sqrt{a}+s}{a^2}  \,
\bigl(\norm{\hat \psi}_\infty + \norm{\triangle^k \hat \psi}_\infty
\bigr)$, where $\triangle = \frac{\partial^2}{\partial \xi_1^2}
+\frac{\partial^2}{\partial \xi_2^2}$ is the frequency domain
Laplacian operator.
\end{prop}

\noindent {\bf Proof.} Observe that, for $t = \left(
\begin{smallmatrix} t_1 \\ t_2
\end{smallmatrix}\right)$ and $x = \left( \begin{smallmatrix} x_1
\\ x_2  \end{smallmatrix}\right)$ in $\R^2$, we have:
\begin{equation*}
\psi_{ast}(x)  =  |\det M|^{-\frac12}  \, \psi( M^{-1}(x-t))
 =  a^{-\frac34} \, \psi \Bigl( \begin{matrix} a^{-1}\, (x_1-t_1)+s \, a^{-1} (x_2-t_2) \\ a^{-\frac12} (x_2-t_2) \end{matrix} \Bigr).
\end{equation*}
The proof then follows from Proposition \ref{pro.loc}, where
\[ R = \{(\xi_1, \xi_2): \xi_1 \in [-\tfrac{2}{a}, -\tfrac{1}{2a}] \cup [\tfrac{1}{2a}, \tfrac{2}{a}],
\, |s - \tfrac{\xi_2}{\xi_1}| \le \sqrt a\}.\] It is easy to check
that $m(R) = \frac{15}{2} \frac{\sqrt{a}+s}{a^2}$. $\qed$

%*****************************************************************************************************

\section{Analysis of singularities} \label{s.ana}

 As observed above, the mother shearlet $\psi$, constructed in
Section~\ref{s.cst}, satisfies $\hat \psi \in C_0^\infty(\widehat
\R^2)$.  It follows that $\psi \in \cS(\R^2)$ and, therefore, the
continuous shearlet transform $ \Sh_f(a,s,t) = \ip{f}{\psi_{ast}}$,
$a>0, s \in \R, t \in \R^2,$ is well defined for all tempered
distributions $f \in \cS'$.

In the following, we will examine the behavior of the continuous
shearlet transform of several distributions containing different
types of singularities. This will be useful to illustrate the basic
properties of the shearlet transform, before stating a more general
result in the next section. Indeed, the rate of decay of the
continuous shearlet transform exactly describes the location and
orientation of the singularities.

In order to state our results, it will be useful to introduce the
following notation to distinguish between the following two
different behaviors of the continuous shearlet transform.

\begin{definition}
Let $f$ be a distribution on $\R^2$, $\Sh_f(a,s,t)$ be defined as in
Definition \ref{defi.shear}, and let $r \in \R$. Then $\Sh_f(a,s,t)$
{\em decays rapidly  as} $a \to 0$, if
\[\Sh_f(a,s,t) = O(a^{k}) \quad \text{ as } a \to 0, \text{ for every } k \ge 0.\]
We use the notation:  $\Sh_f(a,s,t) \sim a^r$ {\em as} $a \to 0$, if
there exist constants $0 < \alpha \le \beta < \infty$ such that
\[\alpha a^r \le \Sh_f(a,s,t) \le \beta a^r \quad \text{ as } a \to 0.\]
\end{definition}

\subsection{Point singularities} \label{s.point}

We start by examining the decay properties of the continuous
shearlet transform of the Dirac $\delta$.

\begin{prop}
If $t = 0$, we have
\[ \Sh_\delta(a,s,t) \sim a^{-\frac34} \quad \mbox{as } a \to 0.\]
In all other cases, $\Sh_\delta(a,s,t)$ decays rapidly as $a \to 0$.
\end{prop}

\smallskip

\noindent {\bf Proof.} For $t=0$ we have
\[ \ip{\delta}{\psi_{ast}} = \psi_{as0}(0) = a^{-\frac34} \, \psi(0) \sim a^{-\frac34}
\; \mbox{ as }a \to 0.\] Next let $t \neq 0$. Then
\[ \ip{\delta}{\psi_{ast}} = \psi_{ast}(0),\]
and, by Proposition~\ref{pro.locshear}, for each $k \in \N$, we have
\[|\psi_{ast}(0)| \le C_k \, a^{-\frac34}(1+a^{-2}t_1^2+2a^{-2}st_1t_2 + (1+a^{-1}
s^2) a^{-1} t_2^2)^{-k}.\]
%&=& C_k a^{-\frac12}(1+a^{-2}t_1^2-2a^{-2}bt_1t_2 + (1+a^{-2}b^2)t_2^2)^{-k}.
Thus, if  $t_2 \ne 0$, then $|\psi_{ast}(0)| =  O(a^{k-3/4})$ as $a
\to 0$. Otherwise, if $t_2 = 0$, $t_1 \ne 0$, then $ |\psi_{ast}(0)|
=  O(a^{2k-3/4})$  as $a \to 0$. $\qed$

%%%%%%%%%%%%%%%%%%%%%%%%%%%%%%%%%%%%%%%%%%%%%%%   OTHER TYPE OF POINT SINGULARITY
%%%%%%%%%%%%%%%%%%%%%%%%%%%%%%%%%%%%%%%%%%%%%%%
\if{

\medskip
Next let us consider the point singularity $\sigma_\alpha(x) =
|x|^{\alpha}$ for $-2 < \alpha < \infty$. The continuous shearlet
transform shows the following decay.

\begin{prop}
Let $\Sh_{\sigma_\alpha}(a,s,t)$ be defined as in Definition
\ref{defi.shear}. If $t = 0$, we have
\[ \Sh_{\sigma_\alpha}(a,s,t) \sim a^{\frac54+\alpha} \quad \mbox{as } a \to 0.\]
In all other cases, $\Sh_{\sigma_\alpha}(a,s,t)$ decays rapidly as
$a \to 0$.
\end{prop}

\smallskip

\noindent {\bf Proof.} First observe that $\hat \sigma_\alpha(x) =
C_\alpha \, |\xi|^{-2-\alpha}$. Using Fubini, we compute
\begin{eqnarray*}
\lefteqn{\ip{\hat{\sigma}_\alpha}{\hat{\psi}_{ast}}}\\
& = & C_\alpha \, \int_{\R^2} |\xi|^{-2-\alpha} \hat \psi_{ast}(\xi) \, d \xi\\
& = & C_\alpha \, a^{\frac34} \int_{\R^2} \hat{\psi}_1(a \xi_1) \,
\hat{\psi}_2(a^{-\frac12}(\tfrac{\xi_2}{\xi_1}-s)) \,
(\xi_1^2+\xi_2^2)^{-1-\frac{\alpha}{2}} e^{-2 \pi i \xi t} \, d \xi_1 \, d \xi_2\\
& = & C_\alpha \, a^{-\frac14} \int_{\R^2} \hat{\psi}_1(\xi_1) \,
\hat{\psi}_2(a^{-\frac12}(\sqrt{a} \,  \tfrac{\xi_2}{\xi_1}-s)) \,
(a^{-2}\xi_1^2+\xi_2^2)^{-1-\frac{\alpha}{2}} \, e^{-2 \pi i (a^{-1}\xi_1,\xi_2)t} \, d \xi_1 \, d \xi_2\\
& = & C_\alpha \, a^{-\frac34} \, \int_{\R^2} \hat{\psi}_1(\xi_1) \,
\hat{\psi}_2(\xi_2) \,
(a^{-2}\xi_1^2+\xi_1^2 a^{-1}(\xi_2+a^{-\frac12} s)^2)^{-1-\alpha/2} \, \\
& &
\cdot e^{-2 \pi i a^{-1/2} (a^{-1/2}\xi_1,\xi_1 (\xi_2+ a^{-1/2} s))t} \, \xi_1 \, d\xi_1 \, d\xi_2\\
& = & C_\alpha \, a^{\frac14+\frac{\alpha}{2}} \int_{\R^2}
\hat{\psi}_1(\xi_1) \, \hat{\psi}_2(\xi_2) \, \xi_1^{-1-\alpha} \,
(a^{-1}+(\xi_2+a^{-\frac12}s)^2)^{-1-\frac{\alpha}{2}} \,\\
&&  \cdot e^{-2 \pi i a^{-1} \xi_1 (t_1+s t_2)} \, e^{-2 \pi i
a^{-1/2} \xi_1 \xi_2 t_2}  d\xi_1 \, d\xi_2 .
\end{eqnarray*}
Let $t = 0$. Then
\[ \ip{\hat \sigma_\alpha}{\hat \psi_{ast}}
= C_\alpha \, a^{\frac14 +\frac{\alpha}{2}} \, \int_{\R}
\frac{\hat{\psi}_1(\xi_1)}{\xi_1^{1+\alpha}} \, \int_\R \frac{\hat
\psi_2 (\xi_2)} {(a^{-1}+(\xi_2+
a^{-\frac12}s)^2)^{1+\frac{\alpha}{2}}} \, d\xi_2 \, d\xi_1.\] For
$a \ll 1$, we obtain
\begin{eqnarray} \label{eq.del}
\left| \int_\R \frac{\hat \psi_2 (\xi_2)} {(a^{-1}+(\xi_2+
a^{-\frac12}s)^2)^{1+\frac{\alpha}{2}}} \, d\xi_2 \right| & = &
a^{1+\frac{\alpha}{2}} \int_\R
\left|\frac{\hat{\psi}_2(\xi_2)}{(1+a(\xi_2+ a^{-\frac12}
s)^2)^{1+\frac{\alpha}{2}}}
\right| \, d\xi_2 \nonumber\\
& \sim & a^{1+\frac{\alpha}{2}} \int_\R |\hat{\psi}_2(\xi_2)| \,
d\xi_2.
\end{eqnarray}
Thus,
\[ \left|\ip{\hat \sigma_\alpha }{\hat \psi_{ast} }\right|
\le C_\alpha \, a^{\frac54+\alpha} \, \int_{\R}
\left|\frac{\hat{\psi}_1(\xi_1)}{\xi_1^{1+\alpha}}\right| \, d\xi_1
\int_\R \left|\hat{\psi}_2(\xi_2)\right| \, d\xi_2 \sim
a^{\frac54+\alpha}  \mbox{ as }a \to 0.\]

Next let $t \neq 0$. If $t_1+ s \,  t_2 \ne 0$, using
again~(\ref{eq.del}), we observe that
\begin{eqnarray*}
\lefteqn{\left|\ip{\hat \sigma_\alpha }{\hat \psi_{ast} }\right|}\\
& \le & C_\alpha \, a^{\frac14+\frac{\alpha}{2}} \, \left| \int_\R
\hat \psi_1 (\xi_1) \, \xi_1^{-1-\alpha} \, e^{-2 \pi i a^{-1} \xi_1
(t_1+s t_2)} \,  \int_\R \frac{\hat \psi_2 (\xi_2)  \, e^{-2 \pi i
a^{-1/2} \xi_1\xi_2 t_2} }
{(a^{-1}+(\xi_2+a^{-\frac12} s)^2)^{1+\frac{\alpha}{2}}} \, d\xi_2 \, d \xi_1 \right| \\
& \le &  C_\alpha' \, a^{\frac54 + \alpha} \, |\widetilde
\psi_1(\tfrac{t_1+st_2}{a})| \, \int_\R
\left|\hat{\psi}_2(\xi_2)\right| \, d\xi_2,
\end{eqnarray*}
where $\widetilde \psi_1(u) =  \int_\R  \xi_1^{-1-\alpha} \,
\hat{\psi}_1(\xi_1)
 \, e^{-2 \pi i u \xi _1} \, d \xi_1$ is a bandlimited function decaying rapidly for $|u| \to \infty$
 (since $\psi_1$ is bandlimited and $C^\infty$, its behavior is similar when $\hat \psi_1(u)$ is divided by $u^{1+\alpha}$).
If $t_1+s \, t_2 = 0$, then one uses a similar estimate employing
the exponential function $e^{-2 \pi i a^{-1/2} \xi_1\xi_2 t_2}$.
$\qed$

Observe that, for $\alpha = -2$, which corresponds to the Dirac
delta, we have the same rate of convergence $a^{-3/4}$ as computed
above.

}\fi
%%%%%%%%%%%%%%%%%%%%%%%%%%%%%%%%%%%%%%%%
%%%%%%%%%%%%%%%%%%%%%%%%%%%%%%%%%%%%%%%%

\subsection{Linear singularities} \label{s.lin}

Next we will consider the linear delta distribution $\nu_p(x_1,x_2)
= \delta(x_1+p \, x_2)$, $p \in \R$, defined by
$$ \ip{\nu_p}{f} = \int_\R f(- p \, x_2,x_2) \, d x_2. $$
The following result shows that the continuous shearlet transform
precisely determines both the position and the orientation of the
linear singularity, in the sense that the transform
$\Sh_{\nu_p}(a,s,t)$ always decays rapidly as $a \to 0$  {\it
except} when $t$ is on the singularity and $s = p$, i.e., the
direction perpendicular to the singularity.

\begin{prop} \label{p.lin}
If $t_1=-p \, t_2 $ and $s=p$, we have
\[ \Sh_{\nu_p}(a,s,t) \sim a^{-\frac14} \quad \mbox{as } a \to 0.\]
In all other cases, $\Sh_{\nu_p}(a,s,t)$ decays rapidly as $a \to
0$.
\end{prop}

\smallskip

\noindent {\bf Proof.} The following heuristic argument gives
\begin{eqnarray*}
\hat \nu_p(\xi_1, \xi_2) & = & \int \int   \delta(x_1+p \, x_2) \, e^{-2 \pi i \xi x} \, d x_2 \, d x_1 \\
& = & \int   e^{-2 \pi i x_2(\xi_2 - p \, \xi_1)} \, d x_2 =
\delta(\xi_2 - p \, \xi_1) = \nu_{(-\frac 1p)} (\xi_1, \xi_2).
\end{eqnarray*}
That is, the Fourier transform of the linear delta on $\R^2$ is
another linear delta on $\widehat{\R}^2$, where the slope $- \frac
1p$ is replaced by the slope $p$. A direct computation gives:
\begin{eqnarray*}
\ip{\hat{\nu}_p \, }{\hat \psi_{ast}}
& = & \int_\R \hat \psi_{ast} (\xi_1, p \xi_1) \, d \xi_1\\
& = & a^{\frac34} \, \int_\R  \hat{\psi}(a \xi_1, \sqrt a  p
\xi_1-\sqrt a s \xi_1) \, e^{-2 \pi i \xi_1(t_1+p t_2)}
\, d \xi_1\\
& = & a^{-\frac14} \,  \int_\R  \hat{\psi}(\xi_1,
a^{-\frac12}p\xi_1-a^{-\frac12} s \xi_1) \, e^{-2 \pi i a^{-1}
\xi_1(t_1+p t_2)} \, d \xi_1\\
& = & a^{-\frac14} \,  \int_\R  \hat{\psi}_1(\xi_1) \,
\hat{\psi}_2(a^{-\frac12}(p-s)) \,
 e^{-2 \pi i a^{-1} \xi_1(t_1+p t_2)} \, d \xi_1\\
& = & a^{-\frac14} \,  \hat{\psi}_2(a^{-\frac12}(p-s)) \,
\psi_1(-a^{-1}(t_1+p t_2) ).
\end{eqnarray*}
If $s \neq p$, then there exists some $a > 0$ such that $|p-s| > \sqrt a$.
This implies that $\hat{\psi}_2(a^{-1/2}(p-s)) = 0,$ and so
$\ip{\hat{\nu}_p \,}{\hat{\psi}_{ast}}=0$. On the other hand, if
$t_1=-p \, t_2 $  and $s=p$, then
$\hat{\psi}_2(a^{-1/2}( p-s)) = \hat{\psi}_2(0) \ne 0$, and
\[ \ip{\hat{\nu}_p \,}{\hat \psi_{ast}}
= a^{-\frac14} \,  \hat{\psi}_2(a^{-\frac12}(p-s)) \,  \psi_1(0)
\sim a^{-\frac14} \; \mbox{ as }a \to 0.\] If $t_1 \ne - p \, t_2 $,
by Proposition~\ref{pro.loc}, we observe that, for all $k \in \N$,
\begin{eqnarray*}
\lefteqn{ \ip{\hat{\nu}_p \,}{\hat \psi_{ast}}}\\
& \le &  a^{-\frac14} \, \hat{\psi}_2(a^{-\frac12}(p-s)) \,  |\psi_1(a^{-1}(t_1+p t_2) )| \\
& \le & C_k \, a^{-\frac14} \, \hat{\psi}_2(a^{-\frac12}(p-s)) \,
(1+a^{-2} (t_1+p t_2)^2)^{-k} =  O(a^{2k-\frac14})  \mbox{ as }a \to
0. \qed
\end{eqnarray*}

\subsection{Polygonal singularities}\label{s.poly}

Here we consider the characteristic function $\chi_V$ of the cone $V
= \{(x_1,x_2) : x_1 \ge 0, \; qx_1 \le x_2 \le px_1\}$, where
$0<q\le p<\infty$. We have the following result.

\begin{prop} \label{prop.decay.cone}
For $t = 0$, if $s =  -\frac1p$ or $s =  -\frac1q$, we have
\[ \Sh_{\chi_V}(a,s,t) \sim a^{\frac34} \quad \mbox{as } a \to 0,\]
and if $s \neq  -\frac1p$ and $s \neq  -\frac1q$, we have
\[ \Sh_{\chi_V}(a,s,t) \sim a^{\frac54} \quad \mbox{as } a \to 0.\]
For $t \neq 0$, if $s =  -\frac1p$ or $s =  -\frac1q$, we have
\[ \Sh_{\chi_V}(a,s,t) \sim a^{\frac34} \quad \mbox{as } a \to 0.\]
In all other cases, $\Sh_{\chi_V}(a,s,t)$ decays rapidly as $a \to
0$.
\end{prop}

The decay of the continuous shearlet transform of $\chi_V$ is
illustrated in Figure~2. As shown in the figure, the decay of
$\Sh_{\chi_V}(a,s,t)$ exactly identifies the location and
orientation of the singularities. It is interesting to notice that
the orientation of the linear singularities can even be detected
considering only the `point singularity' at the origin.

\setlength{\unitlength}{0.13mm}
\begin{figure}
\begin{center}
  \begin{picture}(415,440)(0,0)
    \put(-15,200){\vector(1,0){447}}
    \put(200,-10){\vector(0,1){450}}
  \put(200,200){\line(3,1){215}}
  \put(430,270){$x_2= q \, x_1$}
  \put(200,200){\line(2,3){145}}
  \put(355,415){$x_2= p \, x_1$}
   \put(221,207){\line(-1,3){7}}
   \put(263,221){\line(-1,3){16}}
   \put(305,235){\line(-1,3){27}}
   \put(347,249){\line(-1,3){38}}
   \put(389,263){\line(-1,3){49}}
   \put(200,200){\vector(-1,3){18}}
   \put(200,200){\vector(3,-2){50}}
   \put(200,200){\vector(-1,-1){43}}
    \put(230,130){$\sim a^\frac34$}
  \put(115,270){$\sim a^\frac34$}
  \put(80,125){$\sim a^\frac54$}
  \put(290,337){\vector(-3,2){50}}
  \put(220,390){$\sim a^\frac34$}
  \put(320,240){\vector(1,-3){18}}
  \put(330,160){$\sim a^\frac34$}
  \end{picture}
\end{center}
\caption{Decay properties of the continuous shearlet transform
$\Sh_{\chi_V}(a,s,t)$.}
\end{figure}
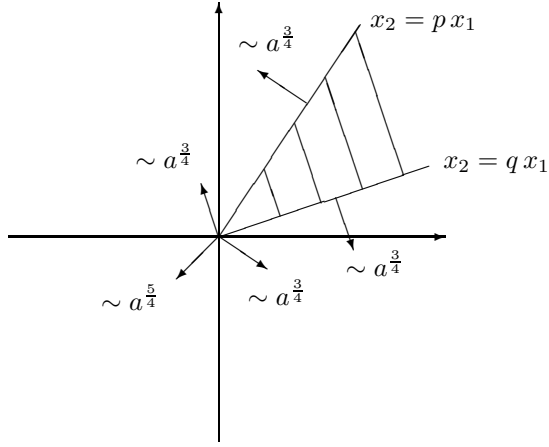

\smallskip

\noindent {\bf Proof of Proposition \ref{prop.decay.cone}.} The
Fourier transform of $\chi_V$ can be computed to be
\[ \hat{\chi}_V(\xi_1,\xi_2) = C\frac{1}{(\xi_1+q\xi_2)(\xi_1+p\xi_2)}, \quad
\mbox{where } C= \frac{(p+q)^2}{(2 \pi)^2}.\]

A direct computation gives:
\begin{eqnarray*}
\lefteqn{\ip{\hat{\chi}_V}{\hat{\psi}_{ast}}}\\ & = & C a^{\frac34}
\int_\R \int_\R \frac{1}{(\xi_1+q\xi_2)(\xi_1+p\xi_2)} \,
\hat{\psi}_1(a\xi_1) \,
\hat{\psi}_2(a^{-\frac12}(\tfrac{\xi_2}{\xi_1}-s)) \, e^{2 \pi i \xi t} \, d\xi_1 \, d\xi_2\\
& = & C a^{-\frac14} \int_\R \int_\R
\frac{1}{(a^{-1}\xi_1+q\xi_2)(a^{-1}\xi_1+p\xi_2)} \,
\hat{\psi}_1(\xi_1)
\, \hat{\psi}_2(a^{-\frac12}(a\tfrac{\xi_2}{\xi_1}-s))\\
& & \cdot e^{2 \pi i (a^{-1}\xi_1,\xi_2) t} \, d\xi_1 \, d\xi_2\\
& = & C a^{-\frac34} \int_\R \int_\R \frac{\xi_1}{(a^{-1}\xi_1+q\xi_1(a^{-1/2}\xi_2+a^{-1}s))(a^{-1}\xi_1+p\xi_1(a^{-1/2}\xi_2+a^{-1}s))} \\
&& \cdot \hat{\psi}_1(\xi_1) \, \hat{\psi}_2(\xi_2) \, e^{2 \pi i \xi_1(a^{-1}(t_1+st_2)+a^{-1/2}\xi_2t_2)} \, d\xi_1 \, d\xi_2\\
& = & C a^{\frac14} \int_\R \int_\R \frac{\xi_1}{(a^{-1/2}\xi_1(1+s q)+q\xi_1\xi_2)(a^{-1/2}\xi_1(1+s p)+p\xi_1\xi_2)}\\
&& \cdot \hat{\psi}_1(\xi_1) \, \hat{\psi}_2(\xi_2)\,  e^{2 \pi i
\xi_1(a^{-1}(t_1+s t_2)+a^{-1/2}\xi_2t_2)} \, d\xi_1 \, d\xi_2.
\end{eqnarray*}
Let us first consider the case $t=0$. By the previous computation we
can rewrite $\ip{\hat{\chi}_V}{\hat{\psi}_{as0}}$ as
\[ C a^{\frac14} \int_\R \int_\R \frac{\xi_1}{(a^{-1/2}\xi_1(1+sq)+q\xi_1\xi_2)(a^{-1/2}\xi_1(1+s p)+p\xi_1\xi_2)}
\, \hat{\psi}_1(\xi_1) \, \hat{\psi}_2(\xi_2) \, d\xi_1 \, d\xi_2.\]
If $s \neq  -\frac1p$ and $s \neq  -\frac1q$, for $a \ll 1$ we can
rewrite $\ip{\hat{\chi}_V}{\hat{\psi}_{as0}}$ as
\begin{eqnarray*}
\lefteqn{\hspace*{-1.9cm}C a^{\frac54} \int_\R \int_\R \frac{\xi_1}{(\xi_1(1+ s
q)+a^{1/2}q\xi_1\xi_2)(\xi_1(1+s p)+a^{1/2}p\xi_1\xi_2)}
\, \hat{\psi}_1(\xi_1) \, \hat{\psi}_2(\xi_2) \, d\xi_1 \, d\xi_2}\\
& \hspace*{-3cm} \sim & \hspace*{-1.5cm}C' a^{\frac54} \int_\R \int_\R \frac{\xi_1}{(\xi_1(1+s
q))(\xi_1(1+s p))} \, \hat{\psi}_1(\xi_1) \, \hat{\psi}_2(\xi_2) \,
d\xi_1 \, d\xi_2,
\end{eqnarray*}
hence
\[ \ip{\hat{\chi}_V}{\hat{\psi}_{as0}} \sim a^{\frac54} \quad \mbox{as }a \to 0.\]
The above computation also shows that if $t=0$ and $s =  -\frac1p$
or $s = -\frac1q$, we have
\[ \ip{\hat{\chi}_V}{\hat{\psi}_{as0}} \sim a^{\frac34} \quad \mbox{as }a \to 0.\]

Next, let us consider the situation, where $t$ lies on one
singularity but $t \neq 0$, i.e., $t_2 = pt_1$ or $t_2 = qt_1$. Here
we will only examine the first case. The second one can be treated
similarly. First let $s = -\frac1p$, i.e., $s$ is perpendicular to
the linear boundary of the cone $x_2 = px_1$. For $a \ll 1$ we have
\begin{eqnarray*}
\ip{\hat{\chi}_V}{\hat{\psi}_{ast}} & = & C a^{\frac14} \int_\R
\int_\R \frac{\xi_1}{(a^{-1/2}\xi_1(1-q/p)+q\xi_1\xi_2)p\xi_1\xi_2}
\hat{\psi}_1(\xi_1) \hat{\psi}_2(\xi_2)\\
&&\cdot e^{2 \pi i a^{-1/2}p t_1 \xi_1 \xi_2} d\xi_1 d\xi_2\\
& = & C a^{\frac34} \int_\R \int_\R
\frac{\xi_1}{(\xi_1(1-q/p)+a^{1/2} q\xi_1\xi_2)p\xi_1\xi_2}\\
&&\cdot\hat{\psi}_1(\xi_1) \hat{\psi}_2(\xi_2)
e^{2 \pi i a^{-1/2}p t_1 \xi_1 \xi_2} d\xi_1 d\xi_2\\
& \sim & a^\frac34\quad \mbox{as }a \to 0.
\end{eqnarray*}
Secondly, let $s \neq -\frac1p$.  We have:
\begin{eqnarray*}
\lefteqn{\ip{\hat{\chi}_V}{\hat{\psi}_{ast}}}\\
& = & C a^{\frac14} \int_\R \int_\R \frac{\xi_1}{(a^{-1/2}\xi_1(1+s q)+q\xi_1\xi_2)(a^{-1/2}\xi_1(1+s p)+a^{-1/2}p\xi_1\xi_2)}\\
&& \cdot  \hat{\psi}_1(\xi_1) \, \hat{\psi}_2(\xi_2) \, e^{2 \pi i \xi_1t_1 (a^{-1}(1+s p) + a^{-1/2}p\xi_2)} \, d\xi_1 \, d\xi_2\\
& = & C a^{\frac14} \int_\R \varphi(\xi_1) \, \hat{\psi}_1(\xi_1) \,
e^{2 \pi i a^{-1}t_1(1+s p)\xi_1} \, d\xi_1,
\end{eqnarray*}
where
\begin{eqnarray*}
\varphi(\xi_1) & = & \int_\R \frac{\xi_1}{(a^{-1/2}\xi_1(1+s q)+q\xi_1\xi_2)(a^{-1/2}\xi_1(1+s p)+p\xi_1\xi_2)}
\hat{\psi}_2(\xi_2)\\
&&\cdot e^{2 \pi i a^{-1/2}t_1p\xi_1\xi_2} d\xi_2.
\end{eqnarray*}
Since $\psi_1$ and $\psi_2$ are band-limited, the function $\varphi$ has
compact support, hence $(\varphi \hat{\psi}_1)^\vee$ is of rapid
decay towards infinity. Thus
\[ \ip{\hat{\chi}_V}{\hat{\psi}_{ast}}
= C \, a^{\frac14} (\varphi \hat{\psi}_1)^\vee(a^{-1}t_1(1-s p)) =
O(a^k) \quad \mbox{as }a \to 0.\] Finally, in case $t_2 \neq pt_1$,
$t_2 \neq qt_1$ and $t_1 \neq 0$, a similar argument to the one
above shows that $\ip{\hat{\chi}_V}{\hat{\psi}_{ast}}$ decays
rapidly  also in this case. $\qed$

\subsection{Curvilinear singularities}\label{s.curv}

We will now examine the behavior of the continuous shearlet
transform of a distribution having a discontinuity along a curve.

Let $B(x_1,x_2) = \chi_{D}(x_1,x_2)$, where $D = \{(x_1,x_2) \in
\R^2 : x_1^2 +x_2^2 \le 1\}$. We have the following:

\begin{prop} \label{prop.decay.disc}
If $t_1^2+t_2^2 =1$ and $s = \tfrac{t_2}{t_1}$, $t_1 \ne 0$, we have
\[ \Sh_{B}(a,s,t) \sim a^{\frac34} \quad \mbox{as } a \to 0.\]
In all other cases, $\Sh_{B}(a,s,t)$ decays rapidly as $a \to 0$.
\end{prop}

The assumption $t_1 \ne 0$ shows that the shearlet transform
$\Sh_{B}(a,s,t)$ is unable to handle the vertical direction $s \to
\infty$. To provide a complete analysis of the singularities of $B$,
we need to use both $\Sh_{B}(a,s,t)$ and $\Sh_{B}^{(v)}(a,s,t)$ (as
defined in Section~\ref{s.cst}). Since the shearlets
$\psi^{(v)}_{ast}$ are defined on the vertical cone $C^{(v)}$, using
$\Sh_{B}^{(v)}(a,s,t)$ one can obtain a similar result to
Proposition \ref{prop.decay.disc}, for $s = \tfrac{t_1}{t_2}$, $t_2
\ne 0$. Since the argument for both cases is exactly the same, we
will only examine the transform $\Sh_{B}(a,s,t)$.

\medskip

In order to prove Proposition~\ref{prop.decay.disc}, we need to
recall the following facts. First, we recall the asymptotic behavior
of Bessel functions, that is given by the following lemma
(cf.~\cite{S93}):

\begin{lemma} \label{l.bessel} There are constants $C_1, C_2$ such that
$$
J_1(2 \pi \lambda) \sim \lambda^{-\frac12} \, (C_1 \, e^{2 \pi i
\lambda} + C_2 \, e^{-2 \pi i \lambda}) \quad \mbox{as }\lambda \to
\infty,$$ and, for $N=1,2, \dots,$ constants $C_1^N, C_2^N$ such
that
$$\Bigl(\frac{d}{d \lambda} \Bigr)^N \, J_1(2 \pi \lambda) \sim \lambda^{-\frac12} \, (C_1^N \, e^{2 \pi i \lambda}
+ C_2^N \, e^{-2 \pi i \lambda}) \quad \mbox{as }\lambda \to
\infty.$$
\end{lemma}
Secondly, we recall the following fact concerning oscillatory
integrals of the First Kind, that can be found in \cite[Ch.8]{S93}:

\begin{lemma} \label{l.curv}
Let $A \in C_0^\infty(\R)$ and $\Phi \in C^1(\R)$, with $\Phi'(t)
\ne 0$ on $\supp A$. Then
$$ I(\lambda) = \int_\R A(t) \, e^{2 \pi i \lambda \Phi(t)} d t = \frac{(-1)^N}{(2 \pi i \lambda)^N} \, \int_\R D^N \bigl( A(t) \bigr) \, e^{2 \pi i \lambda \Phi(t)} d t,$$
for $N=1,2,\dots,$ where $D \bigl( A(t) \bigr)  = \frac d{dt}
\bigl( \frac{A(t)}{\Phi'(t)} \bigr)$.
\end{lemma}

We can now prove Proposition~\ref{prop.decay.disc}.

\smallskip

\noindent {\bf Proof of Proposition \ref{prop.decay.disc}.} The
continuous shearlet transform of $B(x)$ is given by:
\begin{equation}\label{eq.CSTB}
\Sh_{B}(a,s,t) = \ip{B}{\psi_{ast}} = a^{\frac34} \int_\R \int_\R
\hat{\psi}_1(a\xi_1) \,
\hat{\psi}_2(a^{-\frac12}(\tfrac{\xi_2}{\xi_1}-s)) \, e^{-2 \pi i
\xi t} \, \hat B (\xi) \, d\xi_1 \, d\xi_2.
\end{equation}
The Fourier transform $\hat B (\xi_1, \xi_2)$ is the radial
function:
$$ \hat B (\xi_1,\xi_2) = 2 \int_{-1}^1 \sqrt{1-x^2} \, e^{2 \pi i \sqrt{\xi_1^2 + \xi_2^2} x} \, d x =
|\xi|^{-1} \, J_1(2 \pi |\xi|),$$ where $J_1$ is the Bessel function
of order 1. Therefore, the asymptotic behavior of $\hat B (\lambda)$
follows from Lemma \ref{l.bessel}, with the factor $\lambda^{-1/2}$
replaced by $\lambda^{-3/2}$.

Because of the radial symmetry, it is convenient to convert
(\ref{eq.CSTB}) into polar coordinates:
\begin{eqnarray} \label{eq.curv1}
& & \Sh_{B}(a,s,t)     \nonumber \\
 & & \quad = a^{\frac34} \int \! \! \int \hat{\psi}_1(a \rho \cos \theta) \,
\hat{\psi}_2(a^{-\frac12}(\tan \theta-s)) \, e^{-2 \pi i \rho (t_1
\cos \theta + t_2 \sin \theta)} \,
\hat B (\rho) \, \rho \, d \rho \, d \theta \nonumber \\
 & & \quad =  a^{-\frac54} \int \! \! \int \hat{\psi}_1(\rho \cos \theta) \,
\hat{\psi}_2(a^{-\frac12}(\tan \theta-s)) \, e^{-2 \pi i
\frac{\rho}{a} (t_1 \cos \theta + t_2 \sin \theta)} \, \hat B
(\tfrac{\rho}{a}) \, \rho \, d \rho \, d \theta. \,
\end{eqnarray}
We will now examine the asymptotic decay of the function
$\Sh_{B}(a,s,t)$ along the curve $\partial B$
 for $a \to 0$. Thus, we set $t_1^2+t_2^2=1$
 and, without loss of generality,  assume $a<1$.  As we will show, the decay will depend on whether the direction associated with $s$ is normal
 to the curve $\partial B$ or not.

 Let us begin by considering
 the non-normal case $s \ne t_2/t_1$. From (\ref{eq.curv1}), we have:
 \begin{equation*}
   \Sh_{B}(a,s,t)  =
a^{-\frac54} \int I(a,\rho) \, \hat B (\tfrac{\rho}{a}) \, \rho \, d
\rho,
\end{equation*}
where (using the conditions on the support of $\hat \psi_2$)
\begin{equation*}
    I(a,\rho) = \int_{|\tan \theta-s| < \sqrt{a}}  \hat{\psi}_1(\rho \cos \theta) \,
\hat{\psi}_2(a^{-\frac12}(\tan \theta -s)) \, e^{-2 \pi i
\frac{\rho}{a} (t_1 \cos \theta + t_2 \sin \theta)} \, d \theta.
\end{equation*}
Observe that the domain of integration is the cone $| \tan \theta-s|
< \sqrt{a}$ about the direction $\tan \theta = s$, with $a <1$. This
implies that $\theta$ ranges over an interval. Since the conditions
on the support of $\hat \psi_1$ implies that $|\rho \cos \theta|
\subset [\frac12,2]$, it follows that $\rho$ also ranges over an
interval and, as a consequence, $I(a,\rho)$ is compactly supported
in $\rho$.

We will show that $I(a, \rho)$ is an oscillatory integral of the
First Kind that decays rapidly  for $a \to 0$ for each $\rho$. To
show that this is the case, we will apply Lemma~\ref{l.curv} to
$I(a,\rho)$, where $A(\theta;\rho) = \hat{\psi}_1(\rho \cos \theta)
\, \hat{\psi}_2(a^{-1/2}(\tan \theta-s))$, $\Phi(\theta;\rho) = \rho
(t_1 \cos \theta + t_2 \sin \theta)$ and $\lambda = a^{-1}$ and
$\rho$ is a fixed parameter. Observe that $\Phi'(\theta;\rho) = \rho
(-t_1 \sin \theta + t_2 \cos \theta)$ and $\Phi'(\theta;\rho) \ne 0$
for $\tan \theta \ne \frac{t_2}{t_1}$. Thus, we have that, if $|s -
\frac{t_2}{t_1}| \ge \sqrt{a}$, then the function
$\Phi'(\theta;\rho) \ne 0$ on $\supp A$. A direct computation gives
\begin{eqnarray*}
D \bigl( A(\theta;\rho) \bigr) & = & \frac \partial {\partial
\theta} \frac{\hat{\psi}_1(\rho \cos \theta) \,
\hat{\psi}_2(a^{-\frac12}(\tan \theta-s)) }{\rho (-t_1 \sin \theta + t_2 \cos \theta)} =\\
&= & \frac{\sin \theta}{t_1 \sin \theta - t_2 \cos \theta} \hat {\psi_1}' (\rho \cos \theta) \, \hat{\psi}_2(a^{-\frac12}(\tan \theta-s)) + \\
&& + \, a^{-\frac12} \frac{\sec^2 \theta}{\rho(t_2 \cos \theta - t_1 \sin \theta)} \hat {\psi_1} (\rho \cos \theta) \, \hat{\psi}_2'(a^{-\frac12}(\tan \theta-s)) +\\
&& + \, \frac{t_2 \sin \theta + t_1 \cos \theta}{\rho^2(t_2 \cos
\theta - t_1 \sin \theta)^2} \, \hat {\psi_1} (\rho \cos \theta)
\hat{\psi}_2(a^{-\frac12}(\tan \theta-s)).
\end{eqnarray*}
Thus, using the assumptions on $\hat \psi_1, \hat \psi_2$, we have:
$$\bigl| D \bigl( A(\theta;\rho) \bigr) \bigr| < a^{-\frac12} \, C(\theta, \rho) \, \bigl( \norm{\hat{\psi}_1' \hat \psi_2}_\infty +
\norm{\hat{\psi}_1 \hat{\psi}_2'}_\infty + \norm{\hat{\psi}_1 \hat
\psi_2}_\infty \bigr).$$ As observed above, the assumptions on the
support of $\hat \psi_1, \hat \psi_2$ imply that $D \bigl(
A(\theta;\rho) \bigr)$ is compactly supported in $\rho$ away from
$\rho =0$. Using this observation  and $\Phi'(\theta) \ne 0$, we
have that
$$ \norm{ D (A) }_\infty < C \, a^{-\frac12}.$$
Applying the same estimate repeatedly, we have
$$ \norm{ D^N (A) }_\infty < C_N \, a^{-\frac{N}{2}}.$$
Thus, using  Lemma~\ref{l.curv} with $\l = a^{-1}$, we conclude that
$$ \sup_{\rho} |I(a,\rho)| < C \, a^{\frac{N}{2}}. $$
This implies that, under the assumption that we made for $t =
(t_1,t_2)$ and $s$, the function $ \Sh_{B}(a,s,t) $ decays rapidly
for $a \to 0$.

Let us now consider the function $ |\ip{\hat B
}{\hat{\psi}_{ast}}|$, where $t_1^2+t_2^2=1$ and $s= t_2/t_1$
(corresponding to the direction normal to $\partial B$). For
simplicity, let $(t_1,t_2)=(1,0)$. The general case follows using a
similar argument. From (\ref{eq.curv1}), using the change of
variables $u = a^{-1/2} \sin \theta$, we obtain
\begin{equation}\label{eq.curv3}
\ip{\hat B}{\hat{\psi}_{a 0 (1,0)}} =
 a^{-\frac34} \int  \hat B (\tfrac{\rho}{a}) \, \eta_a(\rho) \, \, e^{2 \pi i \frac{\rho}{a}} \rho  \, d \rho
\end{equation}
where
$$
\eta_a(\rho) =  \int_{-(1+a)^{-1/2}}^{(1+a)^{-1/2}}
\hat{\psi}_1(\rho \sqrt{1-au^2}) \,
\hat{\psi}_2(\frac{u}{\sqrt{1-au^2}}) \, e^{2 \pi i \frac{\rho}{a}
(\sqrt{1-au^2} -1)} \, \frac{ d u}{\sqrt{1-au^2}}.$$ The assumptions
on the support of $\hat \psi_2$ imply that
$|\frac{u}{\sqrt{1-au^2}}|<1$. This is equivalent to
$|u|<(1+a)^{1/2}$. Similarly, the assumptions on the support of
$\hat \psi_1$ imply that $|\rho \sqrt{1-au^2}| \subset [\frac12,2]$
and, thus,  $\rho$ ranges over  a closed interval. As a consequence,
the functions $\eta_a(\rho)$ are compactly supported. For $0 < a
<1$, the functions
$$h_a(u) = \hat{\psi}_1(\rho \sqrt{1-au^2}) \,
\hat{\psi}_2(\frac{u}{\sqrt{1-au^2}}) \, \frac{ e^{2 \pi i
\frac{\rho}{a} (\sqrt{1-au^2} -1)}  }{\sqrt{1-au^2}}$$ are
equicontinuous and they converge uniformly:
$$ \lim_{a \to 0} h_a(u) = h_0(u) = \hat{\psi}_1(\rho) \,
\hat{\psi}_2(u) \, e^{-\pi i \rho u^2}.$$ Thus, we have the uniform
limit:
$$ \lim_{a \to 0} \eta_a(\rho) = \eta_0(\rho) = \int_{-1}^{1} \hat{\psi}_1(\rho) \,
\hat{\psi}_2(u) \, e^{-\pi i \rho u^2} \, d \rho,$$ and the same
convergence holds for all $u$-derivatives. In particular,
$\norm{\eta_a}_\infty < C$, for all $a < 1$.

Using the asymptotic estimate given by Lemma~\ref{l.bessel}
into~(\ref{eq.curv3}), we have:
\begin{eqnarray*}
|\ip{\hat B}{\hat{\psi}_{a 0 (1,0)}}| & \sim & a^{-\frac34} \,
\Bigl( C_1 \, \int \Bigl(\frac{a}{\rho}\Bigr)^{\frac32} \eta_a(\rho)
\, e^{4 \pi i \frac{\rho}{a}} \, \rho \, d \rho +  C_2 \, \int
\Bigl(\frac{a}{\rho}\Bigr)^{\frac32}
\eta_a(\rho)  \, \rho \, d \rho \Bigr) \\
&=& a^{\frac34} \, \Bigl( C_1 \, \hat F_a\left(-\frac{2}{a}\right) +
C_2 \int F_a(\rho) \, d \rho \Bigr),
\end{eqnarray*}
where $F_a(\rho) = \eta_a(\rho)  \, \rho^{-1/2}$. The family of
functions $\{F_a: 0 < a< 1\}$ has all its $\rho$ derivatives bounded
uniformly in $a$, and so $\hat F_a(-\frac2a)$ decays rapidly as $a
\to 0$. On the other hand, $\int F_a(\rho) \, d \rho$ tends to $\int
\eta_0(\rho)  \, \rho^{-1/2}\, d \rho$, and thus, combining the two
terms, we obtain that
$$ |\ip{\hat B}{\hat{\psi}_{a 0 (1,0)}}|  \sim a^{\frac34} \quad \mbox{as } a \to 0.$$
If $t$ is not on $\partial B$, then one can show that
$\Sh_{B}(a,s,t) $ has fast decay. This follows from the general
analysis given in Section~\ref{s.wav}. $\qed$

\section{Characterization of the wavefront set using the shearlet transform} \label{s.wav}

The examples described in Section~\ref{s.ana} suggest that the set
of singularities of a distribution on $\R^2$ can be characterized
using the continuous shearlet transform. In this section, we will
show that this is indeed the case. In order to do this, it will be
useful to introduce the notions of singular support and wavefront
set.

For a distribution $u$, we say that $x \in \R^2$ is a {\it regular
point} of $u$ if there is $\phi \in C_0^\infty(U_x)$, where  $U_x$
is a neighborhood of $x$ and $\phi(x) \ne 0$, such that $\phi \, u
\in C_0^\infty(\R^n)$. Recall that the condition $\phi \, u \in
C_0^\infty$ is equivalent to $(\phi \, u)^\wedge$ being rapidly
decreasing. The complement of the regular points of $u$ is called
the {\it singular support} of $u$ and is denoted by $\sing \supp(
u)$. It is easy to see that the singular support of $u$ is  a closed
subset of $\supp( u)$.

The wavefront set of $u$ consists of certain $(x, \l) \in \R^2
\times \R$, with $x \in \sing \supp( u)$. For a distribution $u$, a
point $(x,\l) \in \R^2 \times \R$ is a {\it regular directed point}
for $u$ if there are neighborhoods $U_x$ of $x$ and $V_{\l}$ of
$\l$, and a function $\phi \in C_0^\infty(\R^2)$, with $\phi =1$ on
$U_x$, so that, for each $N>0$ there is  a constant $C_N$ with
$$  |(u \, \phi)^\wedge(\eta)| \le C_N \, (1+|\eta|)^{-N}, $$
for all $\eta = (\eta_1, \eta_2) \in  \R^2$ with
$\frac{\eta_2}{\eta_1}\in V_{\l}.$
 The complement in $\R^2 \times \R$ of the regular directed points for $u$ is called
the {\it wavefront set} of $u$ and is denoted by $WF(u)$. Thus,  the
singular support  is measuring the location of the singularities and
$\l$ is measuring the direction perpendicular to the
singularity.\footnote{This definition is consistent with
\cite{CD04a}, where the direction of the singularity is described by
the angle $\theta$. Observe that our approach does not distinguish
between $\theta$ and $\theta +\pi$, since the continuous shearlets
have frequency support that is symmetric with respect to the origin.
However, in Section \ref{s.ext} we discuss a variant of the
continuous shearlet transform, which can distinguish these cases.}

In the examples presented in Section~\ref{s.ana}, one can verify the
following:
\begin{enumerate}
\item[(i)] { \it Point Singularity $\d(x)$:}  \ $\sing \supp( \delta ) = \{0\}$ and  $WF(\delta) = \{0\} \times \R$.
\item[(ii)] { \it Linear Singularity $\nu_p(x)$:}  $\sing \supp( \nu_p ) = \{(-p x_2, x_2): x_2 \in \R\}$ and
$WF(\nu_p) = \{((-p x_2, x_2), p):\,  x_2 \in \R\}$.
\item[(iii)] { \it Curvilinear Singularity $B(x)$:}   $\sing \supp( B) = \{(x_1,x_2):  \, x_1^2+x_2^2 =1 \}$ and
 $WF(B) = \{((x_1,x_2),\l):  \, x_1^2+x_2^2 =1, \, \l = \tfrac{x_2}{x_1} \}$.
\end{enumerate}

As observed in Section~\ref{s.ana}, all these sets are exactly
identified by the decay properties of the continuous shearlet
transform. Indeed, we have the following general result:

\begin{theorem} \label{th.wavefront}
\begin{enumerate}
\item[(i)] Let $\cR = \{t_0 \in \R^2 :$ for $t$ in a neighborhood
$U$ of $t_0$, $|\Sh_f(a,s,t)| =O(a^k)$ and $|\Sh_f^{(v)}(a,s,t)|
=O(a^k)$ as $a \to 0$, for all $k \in \N$, with the
$O(\cdot)$--terms uniform over $(s,t) \in [-1,1] \times U$\}. Then
\[ \mbox{sing supp}(f)^c = \cR.\]
\item[(ii)] Let $\D = \D_1 \cup \D_2$, where
$\D_1 =\{(t_0,s_0) \in \R^2 \times [-1,1] :$ for $(s,t)$ in a
neighborhood $U$ of $(s_0,t_0)$, $|\Sh_f(a,s,t)| =O(a^k)$ as $a \to
0$, for all $k \in \N$, with the $O(\cdot)$--term uniform over
$(s,t) \in U\}$ and $\D_2=\{(t_0, s_0) \in \R^2 \times [1,\infty) :$
for $(\frac 1s,t)$ in a neighborhood $U$ of $(s_0,t_0)$,
$|\Sh^{(v)}_f(a,s,t)| =O(a^k)$ as $a \to 0$, for all $k \in \N$,
with the $O(\cdot)$--term uniform over $(\frac 1s,t) \in U$\}.
  Then
\[ \mbox{WF}(f)^c = \D.\]
\end{enumerate}
\end{theorem}

The statement (ii) of the theorem shows that the continuous shearlet
transform $\Sh_f(a,s,t)$ identifies the wavefront set for directions
$s$ such that $|s|=|\frac{\xi_2}{\xi_1}| \le 1$ (in the frequency
domain). The continuous shearlet transform $\Sh_f^{(v)}(a,s,t)$
identifies the wavefront set for directions $s$ such that $|s|
=|\frac{\xi_1}{\xi_2}| \le 1$, corresponding to
$|\frac{\xi_2}{\xi_1}| \ge 1$. The proof of
Theorem~\ref{th.wavefront} will require several lemmata and will
adapt several ideas from \cite{CD04a}. The following lemma shows
that if  $t$ is outside  the support of a function $g$, then the
continuous shearlet transform of $g$ decays rapidly as $a \to 0$.

\begin{lemma}
\label{loclemma1} Let $g \in L^2(\R^2)$ with $\norm{g}_\infty  <
\infty$. If $\supp( g ) \subset \cB \subset \R^2$, then for all $k
> 1$,
\[ |\Sh_g(a,s,t)| = \absip{g}{\psi_{ast}} \le C_k\,C(s)^{2}\, \norm{g}_\infty \, a^{\frac14} \, \left(1 + C(s)^{-1}a^{-1}d(t,\cB)^2\right)^{-k},\]
where $C(s) = \left(1 + \frac{s^2}{2} + \left(s^2+\frac{s^4}{4}\right)^\frac12\right)^\frac12$
and $C_k$ is as in Proposition~\ref{pro.locshear}.
\end{lemma}

{\bf Proof.} Since $\norm{g}_\infty < \infty $, by
Proposition~\ref{pro.locshear}, for all $k \in \N$, there is a $C_k
>0$ such that:
\begin{eqnarray}
 \absip{g}{\psi_{ast}}
& \le & \norm{g}_\infty \int_\cB |\psi_{ast}(x)| dx \nonumber \\
& \le &  C_k \, \norm{g}_\infty \, a^{-\frac34} \, \int_\cB  \Bigl(1+ \norm{ M^{-1} (x-t)}^2 \Bigr)^{-k} dx \nonumber\\
& = &  C_k \, \norm{g}_\infty \, a^{-\frac34} \, \int_{\cB+t}
\Bigl(1+ \norm{ M^{-1} x}^2 \Bigr)^{-k} dx, \label{eq.loc1}
\end{eqnarray}
where $M =   \left(\begin{smallmatrix}
     a & -\sqrt a \, s\\ 0 & \sqrt a
 \end{smallmatrix}\right)$.
 Observe that $\norm{x} = \norm{ M M^{-1} x} \le
\norm{M}_{op} \, \norm{M^{-1} x}$, and, thus,
$$\norm{M^{-1} x} \ge \frac 1{\norm{M}_{op}} \, \norm{x}.$$
Since $M =  \left(\begin{smallmatrix}
     1 & -s\\ 0 & 1
 \end{smallmatrix}\right) \,
 \left(\begin{smallmatrix}
     a & 0\\ 0 & \sqrt a
 \end{smallmatrix}\right)$
and $\norm{  \left(\begin{smallmatrix}
     a & 0\\ 0 & \sqrt a
\end{smallmatrix}\right)}_{op} = \sqrt{a}$, then
$$\norm{M^{-1} x} \ge C(s)^{-1} \, a^{-\frac12} \, \norm{x},$$
where $C(s) = \norm{ \left(\begin{smallmatrix}
     1 & -s\\ 0 & 1
 \end{smallmatrix}\right)}_{op}$.
Using these observation in (\ref{eq.loc1}), we have that
\begin{eqnarray*}
\absip{g}{\psi_{ast}}
& \le & C_k \, \norm{g}_\infty \, a^{-\frac34} \, \int_{\cB+t}  \Bigl(1+ C(s)^{-2} a^{-1} \norm{x}^2 \Bigr)^{-k} dx \\
& = & C_k \, C(s)^2 \, \norm{g}_\infty \, a^{\frac14} \,
\int_{C(s)^{-1}a^{-1/2} d(t,\cB)}^\infty  (1+  r^2)^{-k}
r \, dr \\
& = & C_k \,C(s)^{2}\, \norm{g}_\infty \, a^{\frac14} \, \Bigl(1+
C(s)^{-2}a^{-1} d(t,\cB)^2 \Bigr)^{-k}.
\end{eqnarray*}

At last, we compute $C(s)$. We have $ \left(\begin{smallmatrix}
     1 & -s\\ 0 & 1
 \end{smallmatrix}\right) \left(\begin{smallmatrix}
     1 & 0\\ -s & 1
 \end{smallmatrix}\right)
 =  \left(\begin{smallmatrix}
     1+s^2 & -s\\ -s & 1
 \end{smallmatrix}\right)$.
The largest eigenvalue of this matrix is $\lambda_{\max} = 1 +
\frac{s^2}{2} + \left(s^2+\frac{s^4}{4}\right)^\frac12$. Thus we have $C(s)
=\left(1 + \frac{s^2}{2} + \left(s^2+\frac{s^4}{4}\right)^\frac12\right)^\frac12$ for all $s \in
\R$. $\qed$

\if {\bf REMARK:} {\bf [I would suggest to omit this! MODIFY!]}The
dependence on $s$ in Lemma~\ref{loclemma1} cannot be removed. If
this was the case, then, for some subset $\cC$ of $\R^2$ with
positive distance to $0$, we would need that
\[B_s \begin{pmatrix}
     x \\ y
 \end{pmatrix} = \begin{pmatrix}
     1 & -s\\ 0 & 1
 \end{pmatrix}\begin{pmatrix}
     x \\ y
 \end{pmatrix}
=\begin{pmatrix}
     x+sy\\y
 \end{pmatrix}\]
has a positive distance from $0$ for all $(x,y) \in \cC$ independent
of $s$. For this, we only need to consider the `inner' boundary of
$\cC$, which `separates' it from $0$. Let $\epsilon > 0$ and
consider the point $(x_\epsilon, \epsilon) \in \cC$. $S_s$ maps this
point to $(x_\epsilon+s\epsilon, \epsilon)$. Obviously, there exists
an $s \in \R$ with $x_\epsilon+s\epsilon=0$. Since $\epsilon$ is
arbitrary, we can always find a shear $s$, which maps a point of
$\cC$ arbitrarily close to $0$. The application of the parabolic
scaling via $a$ afterwards doesn't change this fact. Thus, in order
to obtain an uniform estimate, it is sufficient to assume that $s
\in S\subset \R$, where $S$ is compact. \fi

\medskip
We can now prove the following inclusions.

\begin{prop}
\label{wavefront1} Let $\mathcal{R}$ and $\D$ be defined as in
Theorem~\ref{th.wavefront}. Then:
\begin{enumerate}
\item[(i)] $  \mbox{sing supp}(f)^c \subseteq \cR. $
\item[(ii)] $ \mbox{WF}(f)^c \subseteq \D. $
\end{enumerate}
\end{prop}

{\bf Proof.} (i) \,  Let $t_0$ be a regular point of $f$. Then there
exists $\phi \in C_0^\infty(\R^2)$ with $\phi(t_0) \equiv 1$ on
$B(t_0,\delta)$, which is the ball centered at $t_0$ with radius
$\delta$, such that $\phi f \in C^\infty(\R^2)$. We will show that
$t_0 \in \cR$. For this, we decompose $\Sh_f(a,s,t)$ as
\begin{equation}\label{loceq1} \Sh_f(a,s,t) = \ip{\psi_{ast}}{\phi f} + \ip{\psi_{ast}}{(1-\phi) f}.
\end{equation}
Observe that we have
\[  \absip{\psi_{ast}}{\phi f}
\le a^{\frac34} \int_{\R^2}
|\hat{\psi}_1(a\xi_1)|\,|\hat{\psi}_2(\tfrac{1}{\sqrt{a}}(\tfrac{\xi_2}{\xi_1}-s))|\,
|\widehat{\phi f}(\xi)| \, d\xi.\] In the following, we will
estimate the above integral for $\xi_1 > 0$. The  case $\xi_1 \le 0$
is similar.

Since $\phi \in C_0^\infty(\R^2)$, for each $k \in \N$ there exists
a constant $C_k$ with $|\widehat{\phi f}(\xi)| \le C_k \,
|\xi|^{-2k}$. Using this fact, together the assumptions on the
support of $\hat{\psi}_{ast}$, for $k > 2$, the first term on the
RHS of \eqref{loceq1} can be estimated as follows:
\begin{eqnarray*}
\lefteqn{a^{\frac34} \int_{\R^+\times \R}
|\hat{\psi}_1(a\xi_1)|\,|\hat{\psi}_2(\tfrac{1}{\sqrt{a}}(\tfrac{\xi_2}{\xi_1}-s))|\,
|\widehat{\phi f}(\xi)| \, d\xi}\\
& \le & C_k \norm{\hat{\psi}}_\infty  a^{\frac34} \int_\frac{1}{2a}^\frac{2}{a} \int_{(s-\sqrt{a})\xi_1}^{(s+\sqrt{a})\xi_1} |\xi|^{-2k}\, d\xi_2 \, d\xi_1\\
& \le & C_k 2^{-k}  \norm{\hat{\psi}}_\infty a^{\frac34} \int_\frac{1}{2a}^\frac{2}{a} \xi_1^{-k} \int_{(s-\sqrt{a})\xi_1}^{(s+\sqrt{a})\xi_1} \xi_2^{-k}\, d\xi_2 \, d\xi_1\\
& = &  \frac {C_k 2^{-k}  \norm{\hat{\psi}}_\infty  a^{\frac34}}{1-k} \left((s+\sqrt{a})^{1-k} - (s-\sqrt{a})^{1-k}\right) \int_\frac{1}{2a}^\frac{2}{a} \xi_1^{1-2k} \, d\xi_1\\
& = & \frac {C_k 2^{-k}  \norm{\hat{\psi}}_\infty
a^{\frac34}}{1-k}\left((\sqrt{a}-s)^{1-k} -
(-\sqrt{a}-s)^{1-k}\right)
\frac{1}{1-2k} \left((\tfrac{2}{a})^{2-2k} -  (\tfrac{1}{2a})^{2-2k}\right)\\
& \le & \frac{C_k 2^{-k} \norm{\hat{\psi}}_\infty
a^{\frac34}}{k(2k-1)} (\sqrt{a}-s)^{1-k}(\tfrac{2}{a})^{2-2k}.
\end{eqnarray*}
Thus, the above quantity  is  a $O(a^k)$ as $a \to 0$, uniformly
over $(t,s) \in B(t_0,\frac{\delta}{2}) \times \R$.

Using Lemma \ref{loclemma1}, we have the following estimate for the
second term on the RHS of \eqref{loceq1}:
\[ \absip{\psi_{ast}}{(1-\phi) f} \le C_k\,C(s)^2\,\norm{(1-\phi) f}_\infty \, a^{\frac14} \,
(1 + C(s)^{-1}a^{-1}d(t,B(t_0,\delta)^c)^2)^{-k},\] where $k \in \N$
is arbitrary. Since $\norm{(1-\phi) f}_\infty < \infty$ and $s$ is
bounded, this yields
\[\absip{\psi_{ast}}{(1-\phi) f} = O(a^k) \quad \mbox{as } a \to 0, \]
 uniformly over $(t,s) \in B(t_0,\frac{\delta}{2})
\times [-1,1]$. A similar estimate hold when $\Sh_f(a,s,t)$ is
replaced by $\Sh^{(v)}_f(a,s,t)$. This proves (i).

(ii) \,  Let $(t_0,s_0)$ be a regular directed point of $f$, with
$s_0 \in [-1,1]$. Then there exists a $\phi \in C_0^\infty(\R^2)$
with $\phi(t_0) \equiv 1$ on a ball $B(t_0,\delta_1)$ such that,
for each $k \in \N$ we have $|\widehat{\phi f}(\xi)| =
O((1+|\xi|)^{-k})$ for all $\xi \in \widehat \R^2$  satisfying
$\frac{\xi_2}{\xi_1} \in B(s_0,\delta_2)$. We will prove that
$(t_0,s_0) \in \D$. For this, we decompose $\Sh_f(a,s,t)$ as in
\eqref{loceq1}. The second term on the RHS of \eqref{loceq1} can
be estimated as in the case (i).  For the first term of
\eqref{loceq1}, we only need to show that $\supp \hat{\psi}_{ast}
\subset \{\xi \in \R^2 : \frac{\xi_2}{\xi_1} \in
B(s_0,\delta_2)\}$ for all $(s,t) \in B(s_0,\delta_2) \times
B(t_0,\delta_1)$, since in this cone $\widehat{\phi f}$ decays
rapidly. As above, we only consider the case $\xi_1 > 0$; the case
$\xi_1 \le 0$ is similar. The support of $\hat{\psi}_{ast}$ in
this half plane is given by
\[ \{(\xi_1,\xi_2) : \xi_1 \in [\tfrac{1}{2a}, \tfrac{2}{a}], \xi_2 \in \xi_1[s-\sqrt{a},s+\sqrt{a}]\}.\]
Let $(s,t) \in B(s_0,\delta_2) \times B(t_0,\delta_1)$. The cone
$\{\xi \in \R^2 : \frac{\xi_2}{\xi_1} \in B(s_0,\delta_2)\}$ is
bounded by the lines $\xi_2 = (s_0-\delta_2) \xi_1$ and $\xi_2 =
(s_0+\delta_2) \xi_1$. Now let $(\xi_1,\xi_2) \in \supp
\hat{\psi}_{ast}$. Then, for $a$ small enough, we have
\[  |\tfrac{\xi_2}{\xi_1} - s_0| \le \sqrt{a} \le \delta_2,\]
and this completes the proof in this case. In case $|s_0| \ge 1$
(this corresponds to $|\frac{\xi_2}{\xi_1}|\le 1$), we proceed
exactly as above, using the transform $\Sh_f^{(v)}(a,s,t)$ rather
than $\Sh_f(a,s,t)$.
 $\qed$

For the converse inclusions we need some additional lemmata. For
simplicity of notation, in the following proofs the symbols $C'$ and
$C_k$ are generic constants and may vary from expression to
expression (in the case of $C_k$, the constant depends on $k$).

\begin{lemma}
\label{loclemma2} Let $S \subset \R$ be a compact set, and $g \in
L^2(\R)$ with $\norm{g}_\infty < \infty$. Suppose that $\supp g
\subset \cB$ for some $\cB \subset \R^2$ and define $(\cB^\eta)^c =
\{x \in \R^2 : d(x,\cB) > \eta\}$. Further define $h \in L^2(\R)$ by
\[ \hat{h}(\xi) = \int_0^\infty \int_{(\cB^\eta)^c} \int_S  \Sh_g(a,s,t) \,
\hat{\psi}_{ast}(\xi) \, ds \, dt \, \frac{da}{a^3}.\] Then
$\hat{h}(\xi)$ decays rapidly as $|\xi| \to \infty$ with constants
dependent only on $\norm{g}_\infty$ and $\eta$.
\end{lemma}

{\bf Proof.}  Using the fact that $S$ is compact, Lemma
\ref{loclemma1} implies that, for each $k > 0$,
\[ |\Sh_g(a,s,t)| \le C_k \, a^{\frac14} \, (1 +  a^{-1}d(t,\cB)^2)^{-k},\]
where $C_k$ depends on $\norm{g}_\infty$ but not on $s$.
 By definition,  the support of
$\hat{\psi}_{ast}$ is contained in the set
\begin{equation}\label{eq.suppG}
\Gamma(a,s) = \{\xi \in \R^2 : \tfrac12 \le a |\xi| \le 2, \:
|s-\tfrac{\xi_2}{\xi_1}| \le \sqrt{a}\}.
\end{equation}
Thus,  $|\hat{\psi}_{ast}(\xi)| \le C' a^{\frac34}
\chi_{\Gamma(a,s)}(\xi)$ and
\begin{equation}\label{loceq2} \int_S \chi_{\Gamma(a,s)}(\xi) \, ds
\le \int_{S \cap \left[\frac{\xi_2}{\xi_1}-\sqrt{a} ,
\frac{\xi_2}{\xi_1}+\sqrt{a}\right]} ds \le C' \, \sqrt{a}.
\end{equation}
Collecting the above arguments,
\begin{eqnarray*}
\hat{h}(\xi) & \le & \int_0^\infty \int_{(\cB^\eta)^c} \int_S
|\Sh_g(a,s,t)| \,
|\hat{\psi}_{ast}(\xi)| \, ds \, dt \, \frac{da}{a^3}\\
& \le & C_k \int_0^\infty \int_{(\cB^\eta)^c} \int_S
a \, \chi_{\Gamma(a,s)}(\xi)\, (1 +  a^{-1}d(t,\cB)^2)^{-k} \, ds \, dt \, \frac{da}{a^3}\\
& \le & C_k\, \int_0^\infty \int_{(\cB^\eta)^c} \int_S
\chi_{\Gamma(a,s)}(\xi)\, ds\,
a \,  (1 +  a^{-1}d(t,\cB)^2)^{-k} \, dt \, \frac{da}{a^3}\\
& \le & C_k\, \int_{\frac{1}{2|\xi|}}^{\frac{2}{|\xi|}} a^{-\frac32}
\, \int_{(\cB^\eta)^c}
(1 +  a^{-1}d(t,\cB)^2)^{-k} \, dt \, da\\
& \le & C_k\, \int_{\frac{1}{2|\xi|}}^{\frac{2}{|\xi|}} a^{-\frac32}
\, \int_\eta^\infty
(1 +  a^{-1}r^2)^{-k} \, r\,dr \, da\\
& \le & C_k \,\int_{\frac{1}{2|\xi|}}^{\frac{2}{|\xi|}} a^{-\frac12}
\,
(1 +  a^{-1}\eta^2)^{-k+2}\, da\\
& \le & C_k \, |\xi|^{-\frac12} \,(1 +  |\xi| \, \eta^2)^{-k+2}.\\
\end{eqnarray*}
Since this holds for each $k > 0$, $\hat{h}(\xi)$ decays rapidly as
$|\xi| \to \infty$.  $\qed$

\begin{lemma}
\label{loclemma3} Let $S \subset \R$ and $\cB \subset \R^2$ be
compact sets. Suppose that $G(a,s,t)$ decays rapidly as $a \to 0$
uniformly for $(s,t) \in S \times \cB$. Define $h \in L^2(\R)$ by
\[ \hat{h}(\xi) = \int_0^\infty \int_\cB \int_S G(a,s,t) \,
\hat{\psi}_{ast}(\xi) \, ds \, dt \, \frac{da}{a^3}.\] Then
$\hat{h}(\xi)$ decays rapidly as $|\xi| \to \infty$.
\end{lemma}

{\bf Proof.} As in Lemma~\ref{loclemma2}, we will use the fact that
$|\hat{\psi}_{ast}(\xi)| \le C' \, a^{\frac34}
\chi_{\Gamma(a,s)}(\xi)$, where $\Gamma(a,s)$ is given by
(\ref{eq.suppG}), and the estimate \eqref{loceq2}. Also, by
hypothesis, for each $k > 0$ and $a >0$ we have
\[ \sup\{|G(a,s,t)| : |\xi| \in [\tfrac{1}{2a},\tfrac{2}{a}], \: t \in \cB\} \le C_k \, a^k.\]
Using all these observation, we have that, for each $k > 0$,
\begin{eqnarray*}
|\hat{h}(\xi)| & \le & \int_0^\infty \int_\cB \int_S |G(a,s,t)| \,
|\hat{\psi}_{ast}(\xi)| \, ds \, dt \, \frac{da}{a^3} \\
& \le & C_k \int_0^\infty \int_\cB \int_S \chi_{\Gamma(a,s)}(\xi)\,a^{k-\frac94}\, ds \, dt \, da\\
& \le & C_k\,  \int_{\frac{1}{2|\xi|}}^{\frac{2}{|\xi|}} a^{k-\frac74}\, da\\
& \le & C_k \, |\xi|^{-k+\frac74}.\qed
\end{eqnarray*}

The proof of the following lemma adapts several ideas from
\cite[Lemma~2.3]{CD05}.

\begin{lemma} \label{loclemma4a}

Suppose $0 \le a_0 \le a_1 < 1$ and $|s| \le s_0$. Then for $K>1$,
there is a constant $C_K$, dependent on $K$ only, such that:
$$ \left| \ip{ \psi_{a_0 s t}}{\psi_{a_1  s' t'}} \right| \le C_K \,
 \left(1 + \frac{a_1}{a_0}\right)^{-K}
  \left(1 + \frac{|s-s'|^2}{a_1}\right)^{-K}
  \left(1 + \frac{\norm{(t-t') }^2}{a_1}\right)^{-K}$$
\end{lemma}

{\bf Proof.} By the properties of $\psi$, we have that, for
$\norm{\xi} > \tfrac 12$ and any $k>0$, there is a corresponding
constant $C_k$ such that
$$ |\hat \psi(\xi)| \le C_k \, \frac 1{(1 + |\xi_1| +|\xi_2|)^{k}}. $$
 We also have that $ \hat \psi(\xi) = 0 $ for
$\norm{\xi} < \tfrac 12$.
 Thus, observing that $M^t_{as} \xi  = (a \, \xi_1,  \sqrt a \, \xi_2-\sqrt a \, s \, \xi_1 )$, it follows that:
$$ |\hat \psi_{ast}(\xi)| \le C_k \, \frac{a^{\frac34}}{(1 + a  \, |\xi_1| + \sqrt a \, |\xi_2-s \, \xi_1 |)^{k}}. $$
Using polar coordinates, by writing $\xi_1 = r \, \cos \theta$ and
$\xi_2 = r \, \sin \theta$, this expression can be written as
$$ |\hat \psi_{ast}(r,\theta)| \le C_k \, \frac {a^{\frac34}}{(1 + a \, r |\cos \theta| + \sqrt a \, r |\sin \theta - s \, \cos \theta |)^{k}}. $$
For $|\theta| \le \pi/2$, using the assumption $|s| \le s_0$, the
last expression can be controlled by
\begin{equation}\label{eq.a1}
|\hat \psi_{ast}(r,\theta)| \le C_k \, \frac {a^{\frac34}}{\left(1 +
a \, r  + \sqrt a \, r \, |\sin \theta-s|\right)^{k}}.
\end{equation}
In addition, since $\sin \theta \sim \theta$ on $|\theta| \le
\pi/2$, we can replace $\sin \theta$ with $\theta$ in the above
expression.
%%%%%%%%%%%%%%%%%%%%%%%%%%%%%%%%%%%%%%%%%%%%%%%%%%%%%%%%%%%%%%
%{\tt {\it Proof of (\ref{eq.a1})}. The claim follows if we show
%that there is a constant $K$ such that:
%$$ K \, (1 + u (|\cos \theta| + a^{-1/2} \, |-s \, \cos \theta + \sin \theta|) \ge 1 + u \left(1 + a^{-1/2} \, |-s+ \sin \theta|\right), $$
%where $u >0$. The last inequality is equivalent to
%$$ K \, |\cos \theta| + a^{-1/2} \, \left( K \, |-s \, \cos \theta + \sin \theta| - |-s +\sin \theta| \right) \ge \frac{1 -K}u +1. $$
%By considering separately the case $|\theta| \le \tfrac{\pi}{4}$
%and $\tfrac{\pi}{2} \le |\theta| < \tfrac{\pi}{4}$, it is now easy
%to verify that the above inequality is satisfied for $K > \sqrt 2
%\, s_0 + 1 +\sqrt 2$. A better bound can be found but it is
%irrelevant for our purpose. }
%%%%%%%%%%%%%%%%%%%%%%%%%%%%%%%%%%%%%%%%%%%%%%%%%%%%%%%%%%%%%%%%

Let $\Delta s = s - s'$, $a_0 = \min(a,a')$ and $a_1 = \max(a,a')$.
Using~(\ref{eq.a1}), and applying the same argument on $|\theta| \le
\pi/2$ and $\pi/2 <|\theta| \le  \pi$, we have that
\begin{eqnarray*} %\label{eq.a2}
 & & \int_{\widehat \R^2} \left|\hat \psi_{a's't'}(\xi) \, \hat \psi_{ast}(\xi) \right| \, d \xi \\
 & & \le C_k \, \int_{\tfrac 1{a_0}}^\infty \int_{-\pi}^{\pi} \frac {(a\, a')^{\frac34} \, r}{\left(1 + a \, r  + \sqrt a \, r \, |\theta-s|\right)^{k}
 \left(1 + a' \, r  + \sqrt{ a'} \, r \, |\theta-s'|\right)^{k}} \, d \theta \, d r \\
 & & \le C_k \,  \int_{a_0^{-1}}^\infty \frac {(a\, a')^{3/4} r}{\left(1 + a \, r\right)^{k'} \left(1 + a' \, r\right)^{k}}
 \int_{-\infty}^\infty \frac {1}{\left(1 + \a \, |\theta|\right)^{k} \left(1 + \a' \, |\theta +\Delta s|\right)^{k}} \, d \theta \, d r,
\end{eqnarray*}
where $\a = \frac{\sqrt a \, r}{1 + a \, r}$, $\a' = \frac{\sqrt{a'}
\, r}{1 + a' \, r}.$ Using a simple calculation we have that, for
$\a > \a' >1$, $k>1$,
\begin{equation*}
\int_{-\infty}^\infty \frac {1}{\left(1 + \a \, |\theta|\right)^{k}
\left(1 + \a' \, |\theta +\Delta s|\right)^{k}} \, d \theta \le C_k
\, \frac 1{\a \, \left(1 + \a' \, |\Delta s|\right)^{k}}.
\end{equation*}
From the definition $\a'$, for $r \ge 1/a'$, we have that  $\tfrac
12 \, \tfrac 1{\sqrt{a'}}\le  \a' \le \tfrac 1{\sqrt{a'}}$. Thus,
for  $r \ge 1/a'$, provided $k>1$, the last expression gives
\begin{equation}\label{eq.a4}
\int_{-\infty}^\infty \frac {1}{\left(1 + \a \, |\theta|\right)^{k}
\left(1 + \a' \, |\theta +\Delta s|\right)^{k}} \, d \theta \le C_k
\,  \frac{(1 + a \, r)}{\sqrt a \, r} \, \left(1 + \frac{|\Delta
s|}{\sqrt{a'}} \right)^{-k}.
\end{equation}
Another simple estimate, provided $k' >1$, gives:
\begin{equation} \label{eq.a5}
\int_{\tfrac 1{a_0}}^\infty \frac {1}{\left(1 + a_0 \, r\right)^{k'}
\left(1 + a_1 \, r\right)^{k}} \, d r \le C_{k'} \, \frac 1{a_0} \,
\left(1 + \frac{a_1}{a_0}\right)^{-k}.
\end{equation}
%%%%%%%%%%%%%%%%%%%%%%%%%%%%%%%%%%%%%%%
%{\tt {\it Proof of (\ref{eq.a5})}. We have that
%\begin{eqnarray*}
%\int_{\tfrac 1{a_0}}^\infty \frac {1}{\left(1 + a_0 \,
%r\right)^{M'} \left(1 + a_1 \, r\right)^{M}} \, d r &\le&
% \left(1 + \frac{a_1}{a_0}\right)^{-M} \int_{\tfrac 1{a_0}}^\infty \frac {1}{\left(1 + a_0 \, r\right)^{M'} } \, d r \\
%& = & \left(1 + \frac{a_1}{a_0}\right)^{-M} \frac 1{a_0} \,
%\int_{1}^\infty \frac {1}{\left(1 + u\right)^{M'} } \, d u.
%\end{eqnarray*}
%Inequality (\ref{eq.a5}) then follows since $M' >1$. }
%%%%%%%%%%%%%%%%%%%%%%%%%%%%%%%%%%%%%%%%%

Thus, using (\ref{eq.a4}) and (\ref{eq.a5}), for $\a > \a' >1$,
$a_0=a'$ and $a_1=a$,  we have that
\begin{eqnarray*}
  \int_{\widehat \R^2} \left|\hat \psi_{a's't'}(\xi) \, \hat \psi_{ast}(\xi) \right| \, d \xi
 &\le& C_k \, \left(\frac{a_1}{a_0}\right)^{\frac14} \, \left(1 + \frac{a_1}{a_0}\right)^{-k+1} \, \left(1 + \frac{|\Delta s|}{\sqrt{a_0}} \right)^{-k} \\
&\le& C_k \,  \left(1 + \frac{a_1}{a_0}\right)^{-k+2} \, \left(1 +
\frac{|\Delta s|}{\sqrt{a_0}} \right)^{-k}.
\end{eqnarray*}
Similarly, for $\a > \a' >1$, $a_1=a'$ and $a_0=a$,  a similar
calculation gives
\begin{eqnarray*}
  \int_{\widehat \R^2} \left|\hat \psi_{a's't'}(\xi) \, \hat \psi_{ast}(\xi) \right| \, d \xi
 &\le& C_k \, \left(\frac{a_1}{a_0}\right)^{\frac34} \, \left(1 + \frac{a_1}{a_0}\right)^{-k} \, \left(1 + \frac{|\Delta s|}{\sqrt{a_1}} \right)^{-k} \\
&\le& C_k \,  \left(1 + \frac{a_1}{a_0}\right)^{-k+1} \, \left(1 +
\frac{|\Delta s|}{\sqrt{a_1}} \right)^{-k}.
\end{eqnarray*}
In general, renaming the index $k$, we can show that
\begin{equation}\label{eq.a6}
  \int_{\widehat \R^2} \left|\hat \psi_{a's't'}(\xi) \, \hat \psi_{ast}(\xi) \right| \, d \xi
 \le C_k \,  \left(1 + \frac{a_1}{a_0}\right)^{-k} \, \left(1 + \frac{|\Delta s|}{\sqrt{a_1}} \right)^{-k} \\
\end{equation}

We have:
$$ \frac{\partial}{\partial \xi_1}  \hat \psi_{ast}(\xi) = (a -\sqrt a \, s ) \, \hat \psi_{ast}(\xi), \quad \frac{\partial}{\partial \xi_2}  \hat \psi_{ast}(\xi) = \sqrt a  \, \hat \psi_{ast}(\xi), $$
$$ \frac{\partial^2}{\partial \xi_1^2}  \hat \psi_{ast}(\xi) = (a -\sqrt a \, s )^2 \, \hat \psi_{ast}(\xi), \quad  \frac{\partial^2}{\partial \xi_2^2}  = a  \, \hat \psi_{ast}(\xi).$$
Thus, observing that $a, a' <1$ and $|s| < s_0$, we have that
$$ \left|\Delta_\xi \, \hat \psi_{ast}(\xi) \, \overline{\hat \psi_{a's't'}(\xi)}\right| \le C' \, a_1 \, |\hat \psi_{ast}(\xi)| \, |\hat \psi_{a's't'}(\xi)|.$$
%%%%%%%%%%%%%%%%%%%%
%{\tt In fact, because $ a<1$, then
%$$ (a + \sqrt a \, s)^2 = (a^2 + 2 a \sqrt a s + s^2 a) \le a_1 + 2 s a_1  + s^2 a_1 < C \, a_1.$$
%}
%%%%%%%%%%%%%%%%%%%
Set
$$ L = I - \frac{\Delta_\xi}{ (2 \pi)^2 \, a_1}. $$
On the one hand, for each $k$, we have
\begin{equation*}
 \left| L^k \left( \hat \psi_{ast} \, \overline{\hat \psi_{a's't'}} \right) (\xi) \right| \le C' \, |\hat \psi_{ast}(\xi)| \, |\hat \psi_{a's't'}(\xi)|.
\end{equation*}
%%%%%%%%%%%%%%%%%%%%
%{\tt Because
%$$ \left| L \left( \hat \psi_{ast} \, \overline{\hat \psi_{a's't'}} \right)\right| \le  \left| \hat \psi_{ast} \, \hat \psi_{a's't'} \right| +
%\left|\frac{\Delta_\xi}{ (2 \pi)^2 \, a_1} \hat \psi_{ast} \,
%\overline{\hat \psi_{a's't'}} \right| \le C' \, |\hat
%\psi_{ast}(\xi)| \, |\hat \psi_{a's't'}(\xi)|.$$ }
%%%%%%%%%%%%%%%%%%%
On the other hand
\begin{equation}\label{eq.a8}
 L^k ( e^{-2 \pi i \xi (t-t')} ) = \left(1 + \frac{\norm{t-t'}^2}{a_1}\right)^k \, e^{-2 \pi i \xi (t-t')}.
\end{equation}
Repeated integrations by parts give
\begin{eqnarray*}
\ip{\psi_{ast}}{\psi_{a's't'}} &=&  \int_{\R^2} \hat \psi_{ast}(\xi) \,  \overline{\hat \psi_{a's't'}(\xi)} \, d \xi \\
&=&  \int_{\R^2} (a \, a')^{3/4} \,  \hat \psi(M^t_{as} \xi ) \,  \overline{\hat \psi(M^t_{a' s'} \xi )} \, e^{- 2 \pi i \xi (t -t')}\, d \xi \\
&=&  \int_{\R^2} L^k \left( (a \, a')^{3/4} \,  \hat \psi(M^t_{as}\xi ) \,  \overline{\hat \psi(M^t_{a' s'} \xi )} \right) \, L^{-k} \left(e^{- 2 \pi i \xi (t -t')} \right)\, d \xi. \\
\end{eqnarray*}
Therefore, from the last expression, using
(\ref{eq.a6})--(\ref{eq.a8}), we have
$$ \left|\ip{\psi_{ast}}{\psi_{a's't'}} \right| \le C_k \, \left(1 + \frac{\norm{t-t'}^2}{a_1}\right)^{-k} \,
\left(1 + \frac{a_1}{a_0}\right)^{-k} \, \left(1 + \frac{|\Delta
s|}{\sqrt{a_1}} \right)^{-k}.$$

The proof is completed recalling that, for $m>0$,  $(1+|x|)^{-2m}
\sim (1+|x|^2)^{-m}$. That is,  there are constants $C_1, C_2 >0$
such that $C_1 \, (1+|x|^2)^{-m} \le (1+|x|)^{-2m} \le C_2 \,
(1+|x|^2)^{-m}$.

\medskip
From Lemma~\ref{loclemma4a}, the following result can be easily
deduced.

\begin{lemma} \label{loclemma4}
Let $\phi_1 \in C^\infty(\R^2)$ supported in $B(0,1)$, and define
$\phi(x) = \phi_1(a_\phi^{-1}(x-t))$.
\begin{enumerate}
\item[(i)] Suppose $0 \le \sqrt{a_0} \le \sqrt{a_1} \le a_\phi< 1$. Then
for $K>0$,
$$ \left| \ip{\phi \, \psi_{a_0  s t}}{\psi_{a_1  s' t'}} \right| \le C_K \,
 \left(1 + \frac{a_1}{a_0}\right)^{-K}
  \left(1 + \frac{|s-s'|^2}{a_1}\right)^{-K}
  \left(1 + \frac{\norm{(t-t') }^2}{a_1}\right)^{-K}.$$
\item[(ii)] Suppose $0 \le \sqrt{a_0} \le a_\phi \le \sqrt{a_1}  < 1$, $a_1
\le a_\phi$. Then for $K>0$,
$$ \left| \ip{\phi \,  \psi_{a_0 s t}}{\psi_{a_1  s' t'}} \right| \le C_K \,
 \left(1 + \frac{a_1}{a_0}\right)^{-K}
  \left(1 + \frac{|s-s'|^2}{a_\phi^2}\right)^{-K}
  \left(1 + \frac{\norm{(t-t') }^2}{a_1}\right)^{-K}.$$
\item[(iii)] Suppose $0 \le \sqrt{a_0} \le a_\phi \le a_1  \le \sqrt{a_1} <
1$. Then for $K>0$,
$$ \left| \ip{\phi \,  \psi_{a_0  s t}}{\psi_{a_1  s' t'}} \right| \le C_K \,
 \left(1 + \frac{a_\phi}{a_0}\right)^{-K}
  \left(1 + \frac{\norm{(t-t') }^2}{a_\phi^2}\right)^{-K}.$$
\end{enumerate}
\end{lemma}

\medskip

Now we now complete the proof of Theorem~\ref{th.wavefront}.

{\bf Proof of Theorem~\ref{th.wavefront}.}

Since one direction was proved by Proposition~\ref{wavefront1}, we
only have to prove the inclusions:
\begin{enumerate}
\item[(i)]  $ \cR \subseteq \mbox{sing supp}(f)^c;$
\item[(ii)] $ \D \subseteq \mbox{WF}(f)^c.$
\end{enumerate}

 First we prove (i). Let $t_0 \in \cR$. Then for all
$t \in B(t_0,\delta)$, a ball centered at $t_0$, we have that
$|\Sh_f(a,s,t)| =O(a^k)$ as $a \to 0$, for all $k \in \N$ with the
$O(\cdot)$--term uniform over $(t,s) \in B(t_0,\delta) \times
[-1,1]$. A similar estimate holds for $\Sh_f^{(v)}(a,s,t)$.

Choose $\phi \in C^\infty(\R^2)$ which is supported in a ball
$B(t_0,\nu)$ with $\nu \ll \delta$ and let $\eta =
\frac{\delta}{2}$. Set $g = \phi f$ and consider the decomposition
\[ \widehat{\phi f}(\xi) = \hat{g}_0(\xi) + \hat{g}_1(\xi) + \hat{g}_2(\xi) + \hat{g}_3(\xi)+\hat{g}_4(\xi),\]
where $\hat{g}_0(\xi) = (\phi \, P(f))^\wedge(\xi)$, where $P(f) =
\int_\R \ip{f}{T_b W} \, T_b W \, d b$,  $W$ is the window
function defined by~(\ref{eq.W}),
   and
\[ \hat{g}_i(\xi) = \chi_{C_1}(\xi) \,  \int_{\cQ_i} \hat{\psi}_{ast}(\xi) \, \Sh_g(a,s,t) \, d\mu(a,s,t),\quad i=1,2,\]
\[ \hat{g}_{i+2}(\xi) = \chi_{C_2}(\xi) \,  \int_{\cQ_i} \hat{\psi}_{ast}(\xi) \, \Sh_g^{(v)}(a,s,t) \, d\mu(a,s,t),\quad i=1,2,\]
where $C_1$, $C_2$ are defined after equation~(\ref{eq.W}),
$d\mu(a,s,t) = \frac{da}{a^3} \, ds \, dt$, $\cQ_1 = [0,1] \times
[-2,2] \times B(t_0,\eta)$ and $\cQ_2 = [0,1] \times [-2,2] \times
B(t_0,\eta)^c$.  The term $\hat{g}_0(\xi)$ decays rapidly as
$|\xi| \to \infty$ since $\phi$ and $P(f)$ are in $C^\infty$. The
term $\hat{g}_2(\xi)$ decays rapidly as $|\xi| \to \infty$ by
Lemma \ref{loclemma2}. In addition, by Lemma \ref{loclemma3},
$\hat{g}_1(\xi)$ decays rapidly as $|\xi| \to \infty$ provided
that $\Sh_g$ decays rapidly as $a \to 0$ uniformly over $(t,s) \in
B(t_0,\eta) \times [-2,2]$. We will consider only the analysis of
the terms $\hat{g}_i$, for $i=1,2$; the cases $i=3,4$ are very
similar.

We claim that $\Sh_g$ indeed decays rapidly as $a \to 0$ uniformly
over $B(t_0,\eta) \times [-2,2]$. In order to prove this, we
decompose $f$ as $f = P(f) + P_{C_1} f + P_{C_2} f$, where
$(P_{C_1} f)^\wedge = \hat f \, \chi_{C_1}$ and $(P_{C_2}
f)^\wedge = \hat f \, \chi_{C_2}.$ It is clear that $\Sh_{\phi \,
P(f)}$ decays rapidly by the smoothness of $\phi$ and $P(f)$.
Next, we examine the term $P_{C_1} f$. The analysis of $P_{C_2} f$
is very similar and will be omitted. We use the decomposition
$P_{C_1} f = f_1+f_2$, where
\[ f_i(x) = \int_{\cQ_i} \psi_{ast}(x) \, \Sh_f(a,s,t) \, d\mu(a,s,t),\quad i=1,2.\]
Let us start by considering the term corresponding to $f_1$. We
have:
$$  \Sh_{\phi f_1}(a,s,t) = \ip{\phi \, f_1}{\psi_{ast} } = \int_{\cQ_1} \ip{\phi \, \psi_{ast}}{ \psi_{a's't'}} \, \Sh_f(a',s',t') \, d\mu(a',s',t'). $$
We will decompose $\cQ_1 = \cQ_{10} \cup \cQ_{11} \cup \cQ_{12}$,
corresponding to $a' > \delta$, $a' \le \delta < \sqrt{a'}$ and
$\sqrt{a'} \le \delta$, respectively. In case $\sqrt{a}, \sqrt{a'}
\le \delta$, by Lemma~\ref{loclemma4} we have that
\begin{equation}\label{eq.est}
\left| \ip{\phi \, \psi_{a s t}}{\psi_{a' s' t'}} \right| \le C_K \,
 \left(1 + \frac{a_1}{a_0}\right)^{-K}
  \left(1 + \frac{\norm{(t-t') }^2}{a_1}\right)^{-K}.
\end{equation}
This implies that, for $m>4$ and $K \ge m-1$
\begin{equation}\label{eq.wf3}
\int_0^\delta  \left(1 + \frac{a_1}{a_0}\right)^{-K} (a')^m
\,\frac{d a'}{(a')^3} \le C_{m,K} \, a^{m-2}, \qquad 0 < a < \delta.
\end{equation}
In fact, for $a' = a_0 \le a =a_1$,
\begin{eqnarray*}
\int_0^a  \left(1 + \frac{a}{a'}\right)^{-K} (a')^m \,\frac{d
a'}{(a')^3} =   a^{m-2} \int_0^1  \left(1 + x\right)^{-K} \, d x =
C_{K} \, a^{m-2}.
\end{eqnarray*}
For $a = a_0 \le a' =a_1$,
\begin{eqnarray*}
\int_a^\delta  \left(1 + \frac{a'}{a}\right)^{-K} (a')^m \,\frac{d a'}{(a')^3} &=&   a^{m-2} \int_1^{\delta/a}  x^{m-3} \left(1 + x\right)^{-K} \, d x \\
&\le & a^{m-2} \int_1^{\infty}  x^{m-3} \left(1 + x\right)^{-K} \, d
x = C_{K,m} \, a^{m-2}.
\end{eqnarray*}
Thus, (\ref{eq.wf3}) follows from the last two estimates. Using
(\ref{eq.wf3}) it follows that
\begin{eqnarray*}
\lefteqn{\int_{\cQ_{12}} \ip{\phi \, \psi_{ast}}{ \psi_{a's't'}} \,
\Sh_f(a',s',t') \, d\mu(a',s',t')}\\
 &\le& C' \int_{-2}^2
\int_{B(t_0,\eta)}
\int_0^\delta \left(1 + \frac{a_1}{a_0}\right)^{-K} (a')^m \,\frac{d a'}{(a')^3} \, d t' \, d s' \\
& \le & C_m \, \, a^{m-2},
\end{eqnarray*}
for all $m >4$. Using the other cases of Lemma~\ref{loclemma4}, one
can show  similar estimates for the integrals over the sets
$\cQ_{10}$ and $\cQ_{11}$. This shows that $  \Sh_{\phi f_1}(a,s,t)$
decays rapidly for $a \to 0$ uniformly over $B(t_0,\eta) \times
[-2,2]$.

Let us consider now the term corresponding to $f_2$:
$$  \Sh_{\phi f_2}(a,s,t) = \ip{\phi \, f_2}{\psi_{ast} } = \int_{\cQ_1} \ip{\phi \, \psi_{ast}}{ \psi_{a's't'}} \, \Sh_f(a',s',t') \, d\mu(a',s',t'). $$
We will decompose $\cQ_2 = \cQ_{21} \cup \cQ_{22}$, corresponding to
$\norm{(t-t')} > \eta$ and $\norm{(t-t')}  \le \eta$, respectively.
Observe that, for $\norm{(t-t')} > \eta$ and $K>1$,
$$ \int_{B(t_0,\eta)^c} \left(1 + \frac{\norm{(t-t') }^2}{a_1} \right)^{-K} d t' \le \int_\eta^\infty \left(1 + \frac{r^2}{a_1} \right)^{-K}
r\, d r \le C' \, a_1 \left(1 + \frac{\eta}{a_1}
\right)^{-K+2}\hspace*{-0.2cm}.  $$ Observe that, on the region $\cQ_{21}$, the
function $\Sh_f(a',s',t')$  is bounded by $C' \, (a')^{3/4}$ since
$f$ is bounded. Thus
\begin{eqnarray*}
& & \int_{\cQ_{21}} \ip{\phi \, \psi_{ast}}{ \psi_{a's't'}} \, \Sh_f(a',s',t') \, d\mu(a',s',t') \\
&& \le C' \int_{-2}^2 \int_0^\delta  \int_{B(t_0,\eta)^c}  \left(1 + \frac{\norm{(t-t') }^2}{a_1} \right)^{-K} \, d t' \left(1 + \frac{a_1}{a_0}\right)^{-K} (a')^{3/4}\,\frac{d a'}{(a')^3}  \, d s' \\
&& \le C'  \int_0^\eta   a_1 \, \left(1 + \frac{\eta}{a_1}
\right)^{-K+2} \left(1 + \frac{a_1}{a_0}\right)^{-K} \frac{d
a'}{(a')^{9/4}},
\end{eqnarray*}
and this is of rapid decay, as $a \to 0$, uniformly over $\cQ_{21}$.
As for the region  $\cQ_{22}$, if $t \in B(t_0,\eta)$ and
$\norm{(t-t')} > \eta$, then $t' \in B(t_0,\delta)$ and thus the
function $\Sh_f$  decays rapidly, for $a \to 0$, over this region.
Repeating the analysis as in the case $\cQ_{12}$, we can prove that
$ \int_{\cQ_{22}} \ip{\phi \, \psi_{ast}}{ \psi_{a's't'}} \,
\Psi_f(a',s',t') \, d\mu(a',s',t')$ is of rapid decay, as $a \to 0$,
uniformly over $\cQ_{22}$. Combining these observations, we conclude
that $  \Sh_{\phi f_2}(a,s,t)$ decays rapidly as $a \to 0$ uniformly
over $B(t_0,\eta) \times [-2,2]$.

It follows that  $\Sh_g(a,s,t)$ decays rapidly  as $a \to 0$
uniformly for all $(s,t) \in B(t_0,\eta) \times [-2,2]$ and, thus,
by Lemma \ref{loclemma3}, $\hat{g}_1(\xi)$  decays rapid  as $|\xi|
\to \infty$. We can now conclude that $\hat g$ decays rapidly  as
$|\xi| \to \infty$, hence completing the proof of (i).

For part (ii), we only sketch the idea of the proof, since it is
very similar to part (i). Let $(t_0,s_0) \in \D$. We consider
separately the case $|s_0| \le 1$ and $|s_0| \ge 1$. In the first
case, for all $t \in B(t_0,\delta)$ and $s \in B(s_0, \Delta)$, we
have that $|\Sh_f(a,s,t)| =O(a^k)$, as $a \to 0$, for all $k \in \N$
with the $O(\cdot)$--term uniform over $(t,s) \in B(t_0,\delta)
\times B(s_0, \Delta)$. Choose $\phi \in L^2(\R^2)$ which is
supported in a ball $B(t_0,\nu)$ with $\nu \ll \delta$ and let $\eta
= \frac{\delta}{2}$. Then the proof proceeds as in part (i),
replacing $B(t_0,\delta) \times  [-2,2]$ with $B(t_0,\delta) \times
B(s_0,\Delta)$. Also, for the estimates involving inner products of
$\psi_{ast}$ and $\psi_{a's't'}$ we will now use
Lemma~\ref{loclemma4} including the directionally sensitive term.
For example, when $\sqrt{a}, \sqrt{a'} \le \delta$, by
Lemma~\ref{loclemma4} we will use the estimate
$$ \left| \ip{\phi \, \psi_{a s t}}{\psi_{a' s' t'}} \right| \le C_K \,
 \left(1 + \frac{a_1}{a_0}\right)^{-K} \left(1 + \frac{|s-s'|^2}{a_1}\right)^{-K}
  \left(1 + \frac{\norm{(t-t') }^2}{a_1}\right)^{-K},$$
rather than (\ref{eq.est}). We can proceed similarly for the other
estimates. The proof for the case $|s_0| \ge 1$ is exactly the same,
with the transform $\Sh_f^{(v)}(a,s,t)$ replacing $\Sh_f(a,s,t)$.
$\qed$

%*************************************************************************************************
%*************************************************************************************************
%****************************************************************************************************

\section{Extensions and generalizations of the continuous shearlet transform}
\label{s.ext}

As mentioned above, there are several variants and generalizations
of the continuous shearlets introduced in Section~\ref{s.cst}. In
fact, using  the theory of the affine systems, we can obtain several
other examples of continuous wavelets depending on the three
variables: scale, shear and location.

For example, we can generalize our construction by considering the
case where $\L$ is given by (\ref{def.G}) and  $M$ is of the form
$$ M_\d =\begin{pmatrix}  a & - a^{\delta} \, s \\
\\ 0  & a^{\delta}  \end{pmatrix} = B \, A_\d, \quad   a \in I, s \in S,
$$
where
$ B =  \left(\begin{smallmatrix}  1 & -s \\
\\ 0  & 1  \end{smallmatrix}\right)$,  $ A_\d = \left(\begin{smallmatrix}  a & 0 \\
\\ 0  & a^\d  \end{smallmatrix}\right)$
and  $0 \le \delta \le 1$ is fixed. If  $\delta=\frac 12$, then  we
have the case of continuous shearlets. In general, for other choices
of $\d$, $A_\d$ will provide different kinds of anisotropic scaling.
Using a construction similar to the continuous shearlets
 in Section~\ref{s.cst}, for each $0 \le \delta \le 1$, one can construct systems of the form
 $$ \{\psi_{ast} = T_t \, D_{M_\d} \,\psi: \, (M_\d,t) \in \L\},$$
where $\hat \psi_{ast}$ is supported on the set:
$$ \supp \hat{\psi}_{ast} \subset \{(\xi_1, \xi_2):
\xi_1 \in [-\tfrac{2}{a}, -\tfrac{1}{2a}] \cup [\tfrac{1}{2a},
\tfrac{2}{a}], \, |\tfrac{\xi_2}{\xi_1}-s| \le a^{1-\d}\}.$$ This
provides a new family of transforms $\Sh^\d_f(a,s,t) =
\ip{f}{\psi_{ast}}$, for various values of $\d$. It turns out that,
provided $0 \le \d <1$, the behavior of the transforms
$\Sh^\d_f(a,s,t)$ is very similar to the continuous shearlet
transform in dealing with pointwise and linear singularities. More
precisely, one can repeat the analysis in Sections~\ref{s.point},
\ref{s.lin} and \ref{s.poly} using $\Sh^\d_f(a,s,t)$ (for $\d \ne
1$) rather than the continuous shearlet transform. Similarly, the
behavior of the transforms $\Sh^\d_f(a,s,t)$ is very similar to the
continuous shearlet transform in dealing with curvilinear
singularities (see Section~\ref{s.curv}), provided $0 < \d <1$.
However, the proof of Theorem~\ref{th.wavefront} requires $\d =
\frac 12$; that is, we need to use the continuous shearlet
transform.

Another variant of affine systems generated by $\L$, given by
(\ref{def.G}), is obtained by reversing the order of the shear and
dilation matrices, namely, by letting
$$ M_\d =\begin{pmatrix}  a & - a \, s \\
\\ 0  & a^{\delta}  \end{pmatrix} =A_\d \, B , \quad   a \in I, s \in S,
$$
where $A_\delta$, $B$ are defined as above, and $0 \le \delta \le 1$
is fixed. Also in this case, we can construct variants of the
continuous shearlets. Using the same ideas as above, we obtain a
system of  functions $\psi_{ast}$ with support
$$ \supp \hat{\psi}_{ast} \subset \{(\xi_1, \xi_2):
\xi_1 \in [-\tfrac{2}{a}, -\tfrac{1}{2a}] \cup [\tfrac{1}{2a},
\tfrac{2}{a}], \, |a^{1-\d} \tfrac{\xi_2}{\xi_1}-s| \le
a^{1-\d}\}.$$ It turns out (as one can see from the support
condition) that the transform associated with these systems is not
even able to `locate' the linear singularities, in the sense
described in Section~\ref{s.lin}.

\medskip

As we mentioned above, because the frequency support of the
continuous shearlets is symmetric with respect to the origin, the
continuous shearlet transform is unable to distinguish the
orientation associated with the angle $\theta$ from the angle
$\theta + \pi$. In order to be able to distinguish these two
directions, we can modify our construction as follows. Let $\psi \in
L^2(\R^2)$ be defined as in Section~\ref{s.cst}, except that $\supp
\hat \psi_1 \subset [\frac12,2]$. That is, we have a one-sided
version of the shearlet $\hat \psi$ illustrated in Figure~1. It is
then clear that, if we consider the affine system generated by
$\psi$ under the action of $\L$ given by (\ref{def.G}) and
(\ref{def.G2}), this cannot provide a reproducing system for all of
$L^2(\R^2)$. In fact, $\psi$ is a wavelet for the subspace
$L^2(H)^\vee$, where $H = \{(\xi_1,\xi_2) \in \widehat \R^2: \, \xi_1
\ge 0\}$. In order to obtain a wavelet for the space $L^2(\R^2)$, we
can modify the set $\L$ as follows. Let $\L'$ given by
(\ref{def.G}), where $G \subset GL_2(\R)$ is the set of matrices:
\begin{equation}\label{def.G3}
G = \left\{ M =  \begin{pmatrix} \ell a & - \ell \sqrt a \, s \\
\\ 0  & \ell \sqrt a  \end{pmatrix}, \quad a \in I, \, s \in S, \, \ell=-1,1 \right\},
\end{equation}
where $I \subset \R^+$, $S \subset \R$. Then the function $\psi$ is
a continuous wavelet for $L^2(\R^2)$ with respect to $\L'$. These
modified versions of the continuous shearlets depend on four
variables $a,s,t,\ell$ and have frequency support:
$$ \supp \hat{\psi}_{ast\ell} \subset \begin{cases}   \{(\xi_1, \xi_2):
\xi_1 \in [-\tfrac{2}{a}, -\tfrac{1}{2a}],  \,
|\tfrac{\xi_2}{\xi_1}-s|
\le \sqrt a\}, \, & \text{ if } \ell = -1, \\
  \{(\xi_1, \xi_2):
\xi_1 \in  [\tfrac{1}{2a}, \tfrac{2}{a}], \,
|\tfrac{\xi_2}{\xi_1}-s| \le \sqrt a\}, \, & \text{ if } \ell = 1.
\end{cases}
$$
We remark that these modified versions of the continuous shearlets
are in fact complex functions, whereas the shearlets we use
throughout this paper are real functions.
\medskip

Finally, there exists a natural way to construct continuous
shearlets also in dimensions larger than 2. We refer to
\cite{GLLWW3} for a discussion about the generalizations of
shear matrices to higher dimensions.

\appendix
\section{Additional Computations} \label{app}

{\bf Proof of Theorem~\ref{th.rf}.}

Suppose that (\ref{eq.admis}) holds. Then, by applying Parseval and Plancherel
formulas, for any $f \in L^2(\R^n)$ we have:
\begin{eqnarray*}
 & &\int_{\R^n} \int_{G} \left|\ip{f}{T_t \, D_M \,  \psi}\right|^2 \,   d \l(M) \, d t \\
 & & = \int_{\R^n} \int_{G} \left| \int_{\R^n} \hat f (\xi) \, \overline{ \hat \psi(M^t \xi)} \, e^{2 \pi i \xi t} d \xi \right|^2 \, |\det M| \,   d \l(M) \, d t\\
& & = \int_{\R^n} \int_{G} \left| \left( \hat f \, \overline{ \hat \psi(M^t \cdot)} \right)^\vee(t) \right|^2 \, |\det M| \,   d \l(M) \, d t\\
& & = \int_{G} \int_{\R^n}  \left| \left( \hat f \, \overline{ \hat \psi(M^t \cdot)} \right)^\vee(t) \right|^2 \, d t \, |\det M| \,   d \l(M) \\
& & =  \int_{G} \int_{\R^n} | \hat f(\xi)|^2  \, |\hat \psi(M^t \xi)|^2  \, |\det M| \, d \xi  \,  d \l(M) \\
& & =  \int_{G} | \hat f(\xi)|^2  \, \Delta(\psi)(\xi) \, d \xi = \norm{f}^2.
\end{eqnarray*}
Equation~(\ref{eq.calderon}) follows from the above equality by polarization.

Conversely, suppose that
$$ \int_{\R^n} \int_{G} \left|\ip{f}{T_t \, D_M \,  \psi}\right|^2 \,   d \l(M) \, d t = \norm{f}^2$$
for all $f \in L^2(\R^n)$. Let $\xi_0$ be a point of differentiability of $\Delta(\psi)(\xi)$
and let $\hat f (\xi) = |B(\xi_0,r)|^{-1/2} \, \chi_{B(\xi_0,r)}(\xi)$, where $B(\xi_0,r)$ is a ball centered at $\xi_0$ of radius $r$.
By reversing the chain of equalities above we conclude that
$$ \frac 1{|B(\xi_0,r)|} \,  \int_{B(\xi_0,r)} \Delta(\psi)(\xi) \, d \xi = 1 $$
for all $r>0$. Letting $r \to 0$, we obtain that $\Delta(\psi)(\xi_0) =1$, and, since a.e. $\xi \in \R^n$
is a point of differentiability, (\ref{eq.admis}) holds. $\qed$

\smallskip

The proof easily extends to function $f \in L^2(V)^\vee$. In fact, it suffices to replace $\hat f$ with
$\hat f \chi_V$ in the argument above.

\medskip

{\bf Proof of equality~(\ref{eq.dec3}).}

Using Plancherel and Parseval formulas we have that
\begin{eqnarray}
 \int_{\R} \left|\ip{f}{T_t \, W}\right|^2 \, \, d t  &=& \int_{\R} \left|\int_\R \hat f(\xi) \, \overline{ \widehat W(\xi)} \, e^{2 \pi i \xi t} \, d \xi \right|^2 \, \, d t \nonumber\\
&=& \int_{\R} \left| \left(\hat f \, \overline{ \widehat W} \right)^\vee(t) \right|^2 \, \, d t \nonumber\\
&=& \int_{\R} |\hat f(\xi)|^2  \, |\widehat W(\xi)|^2 \, \, d \xi. \label{eq.A1}
\end{eqnarray}
Using a similar computation we have:
\begin{eqnarray}
& & \int_{\R} \int_{-2}^2 \int_0^1 \left|\ip{P_{C_1} f}{\psi_{ast}}\right|^2 \,   \frac{d a}{a^3} \, ds \,  d t \nonumber\\
& & =
\int_{\R} \int_{-2}^2 \int_0^1 \left| \int_\R \hat f(\xi)  \, \chi_{C_1}(\xi)  \, \overline{\hat \psi(M_{as}^t \xi)} \, e^{2 \pi i \xi t} \, d \xi \right|^2 \,   \frac{d a}{a^{3/2}} \, ds \,  d t \nonumber\\
& & =
\int_{-2}^2 \int_0^1 \int_{\R} \left| \left( \hat f \, \chi_{C_1} \, \overline{\hat \psi(M_{as}^t \cdot)}  \right)^\vee(t) \right|^2 \,  d t \,   \frac{d a}{a^{3/2}} \, ds  \nonumber\\
& & =
\int_{-2}^2 \int_0^1 \int_{\R}  |\hat f(\xi)|^2 \, \chi_{C_1}(\xi) \, |\hat \psi(M_{as}^t \xi)|^2  \,  d \xi \,   \frac{d a}{a^{3/2}} \, ds. \label{eq.A2}
 \end{eqnarray}
As in the proof of Proposition~\ref{prop2}, for $\xi \in C_1$, we have that
$$ \int_{-2}^2 \int_0^1    |\hat \psi(M_{as}^t \xi)|^2  \,  \frac{d a}{a^{3/2}} \, ds = \int_0^1   |\hat \psi_1(a \xi_1)|^2  \,  \frac{d a}{a}.$$
Thus, using the last equality, from (\ref{eq.A2}) we have that
\begin{equation}\label{eq.A3}
\int_{\R} \int_{-2}^2 \int_0^1 \left|\ip{P_{C_1} f}{\psi_{ast}}\right|^2 \,   \frac{d a}{a^3} \, ds \,  d t =
\int_\R  |\hat f(\xi)|^2   \chi_{C_1}(\xi) \int_0^1   |\hat \psi_1(a \xi_1)|^2  \,  \frac{d a}{a} \, d \xi.
\end{equation}
Similarly, we have that
\begin{equation}\label{eq.A4}
\int_{\R} \int_{-2}^2 \int_0^1 \left|\ip{P_{C_2} f}{\psi_{ast}^{(v)}}\right|^2 \,   \frac{d a}{a^3} \, ds \,  d t =
\int_\R  |\hat f(\xi)|^2   \chi_{C_2}(\xi) \int_0^1   |\hat \psi_1(a \xi_2)|^2  \,  \frac{d a}{a} \, d \xi.
\end{equation}
Thus, combining (\ref{eq.A1}), (\ref{eq.A3}) and (\ref{eq.A4}) and using~(\ref{eq.W}) we have that
\begin{eqnarray*}
 && \int_{\R} \left|\ip{f}{T_t \, W}\right|^2 \, \, d t + \int_{\R} \int_{-2}^2 \int_0^1 \left|\ip{P_{C_1} f}{\psi_{ast}}\right|^2 \,   \frac{d a}{a^3} \, ds \,  d t \\
 && + \int_{\R} \int_{-2}^2 \int_0^1 \left|\ip{P_{C_2} f}{\psi_{ast}^{(v)}}\right|^2 \,   \frac{d a}{a^3} \, ds \,  d t \\
 && = \int_\R  |\hat f(\xi)|^2   \left(|\widehat W(\xi)|^2 + \chi_{C_1}(\xi) \int_0^1   |\hat \psi_1(a \xi_1)|^2  \,  \frac{d a}{a}
 +\chi_{C_2}(\xi) \int_0^1   |\hat \psi_1(a \xi_2)|^2  \,  \frac{d a}{a} \right) \, d \xi \\
 && = \norm{f}^2. \qed
 \end{eqnarray*}

%*************************************************************************************************
%*************************************************************************************************
%****************************************************************************************

\section*{Acknowledgments}
The research for this paper was performed while the first author
was visiting the Departments of Mathematics at Washington
University in St.~Louis and the Georgia Institute of Technology.
This author thanks these departments for their hospitality and
support during these visits. We are also indebted to F. Bornemann,
E.~Cand\`{e}s, L.~Demanet, K.~Guo,  W.~Lim,  G.~Weiss, and
E.~Wilson for helpful discussions.

\bibliographystyle{amsplain}

\begin{thebibliography}{10}

\bibitem{BI76} J. Bros and D. Iagolnitzer, {\em Support essentiel et structure analytique des distributions},
S´eminaire Goulaouic-Lions-Schwartz, exp. no. 19 (1975-1976).
\bibitem{Cal64} A. Calder\`on, {\em Intermediate spaces and interpolation. The complex method}, Stud. Math.
\textbf{24} (1964), 113--190.
\bibitem{CD05} E. J. Cand\`es,  and L. Demanet,
{\em The curvelet representation of wave propagators is optimally sparse},
to appear in Comm. Pure Appl. Math. (2005).
\bibitem{CD99} E. J. Cand\`es and D. L. Donoho, {\em Ridgelets: a key to higher--dimensional intermittency?},
       Phil. Trans. Royal Soc. London A \textbf{357} (1999), 2495--2509.
\bibitem{CD04} E. J.~Cand\`es and D. L.~Donoho, {\it New tight frames of curvelets and optimal representations
               of objects with $C^2$ singularities}, Comm. Pure Appl. Math. \textbf{56} (2004), 219--266.
\bibitem{CD04a} E. J.~Cand\`es and D. L.~Donoho, {\it Continuous curvelet transform: I. Resolution of the wavefront set},
Appl. Comput. Harmon. Anal. {\bf 19} (2005), 162--197.
\bibitem{CD04b} E. J.~Cand\`es and D. L.~Donoho, {\it Continuous curvelet transform: II. Discretization of frames},
Appl. Comput. Harmon. Anal. {\bf 19} (2005), 198--222.
\bibitem{Cas00} P. G. Casazza, {\em The art of frame theory}, Taiwanese J. Math., \textbf{4} (2000), 129--201.
\bibitem{Chr03} O.~Christensen,
\emph{An Introduction to Frames and Riesz Bases}, Birkh\"auser,
Boston, 2003.
\bibitem{CF78} A. C\'ordoba,  and C.~Fefferman,
  {\em Wave packets and Fourier integral operators}, Comm.
  Partial Diff. Eq. \textbf{3} (1978), 979--1005.
%\bibitem{Dau92}   I.  Daubechies,     \emph{Ten Lectures on Wavelets}, CBMS--NSF Regional     Conference Series in Applied mathematics, Vol.61, SIAM,     Philadelphia, 1992.
%\bibitem{ELL} G. Easley, D. Labate, and W-Q. Lim
%     {\em Sparse Directional Image Representations using Shearlets}, preprint (2006).
\bibitem{GMP85} A. Grossmann, J. Morlet, and T. Paul, {\em Transforms associated to square integrable
group representations I: General Results}, J. Math. Phys.
\textbf{26} (1985), 2473--2479.
\bibitem{GKL05} K. Guo, G. Kutyniok, and D. Labate,
{\em Sparse Multidimensional Representations using Anisotropic
Dilation and Shear Operators,} in: Wavelets and Splines, G. Chen and
M. Lai (eds.), Nashboro Press, Nashville, TN (2006), 189--201.
\bibitem{GL} K. Guo and D. Labate,
     {\em Optimally Sparse Multidimensional Representation using Shearlets}, preprint (2006).
\bibitem{GLLWW} K. Guo, W-Q. Lim, D. Labate, G. Weiss and E. Wilson,
     {\em Wavelets with composite dilations}, {Electr. Res. Announc. of AMS} \textbf{10} (2004), 78--87.
\bibitem{GLLWW3} K. Guo, W-Q. Lim, D. Labate, G. Weiss and E. Wilson,
     {\em Wavelets with composite dilations and their MRA properties},
     \textsl{Appl. Comput. Harmon. Anal.} \textbf{20} (2006), 220--236.
%\bibitem{HLW02}     E.~Hern\'andez, D.~Labate, and  G.~Weiss,
%      \emph{A unified characterization of reproducing systems
%     generated by a finite family, II},  J. Geom. Anal.
%     \textbf{12(4)} (2002), 615--662.
%\bibitem{HW96}  E.~Hern\'andez and G.~Weiss, \emph{A First Course on Wavelets},     CRC Press, Boca Raton FL, 1996.
\bibitem{Hol95} M. Holschneider, {\em Wavelets. Analysis tool}, Oxford University Press, Oxford, 1995.
\bibitem{Hor} L.~H\"ormander,  \emph{The analysis of linear partial differential operators. I. Distribution theory and Fourier analysis.}
  Springer-Verlag, Berlin, 2003.
\bibitem{LLKW05}  D. Labate, W.-Q. Lim, G. Kutyniok, and G. Weiss,
{\em Sparse multidimensional representation using shearlets},
Wavelets XI (San Diego, CA, 2005), 254--262, SPIE Proc. {\bf 5914},
SPIE, Bellingham, WA, 2005.
\bibitem{LWWW}     R.~S.~Laugesen, N.~Weaver, G.~Weiss, and E.~Wilson, \emph{A
     characterization of the higher dimensional groups associated
     with continuous wavelets}, J. Geom. Anal. \textbf{12(1)} (2001), 89--102.
%\bibitem{LM05} E. Le Pennec and S. Mallat, Sparse geometric image representations with bandelets,
%IEEE Trans. Image Process. \textbf{14} (2005), 423--438.
\bibitem{Mal98} S. Mallat,  {\it A Wavelet Tour of Signal Processing},
    Academic Press, San Diego, 1998.
\bibitem{Mey}  Y.  Meyer, {\it Wavelets and Operators}, Cambridge Stud. Adv. Math. vol. 37, Cambridge Univ. Press,
Cambridge, UK, 1992.
\bibitem{Smi98} H. F. Smith, \emph{A Hardy space for Fourier integral operators},
J. Geom. Anal. {\bf 8}, 629--653.
\bibitem{Sog} C.~D.~Sogge, {\it Fourier Integrals in Classical Analysis,}
Cambridge University Press, Cambridge, 1993.
\bibitem{S93} E. M. Stein,  \emph{Harmonic Analysis: real-variable methods, orthogonality, and oscillatory integrals},
Princeton University Press, Princeton, NJ, 1993.
%\bibitem{SW70} E. M. Stein and G. Weiss, \emph{Introduction to Fourier Analysis on Euclidean Spaces}, Princeton University Press, Princeton, NJ, 1970.
\bibitem{WW} G.~Weiss and E.~Wilson, \emph{The mathematical theory of wavelets},
     Proceeding of the NATO--ASI Meeting. Harmonic Analysis 2000 -- A Celebration.
     Kluwer Publisher, 2001.

\end{thebibliography}

\end{document}